\documentstyle[12pt]{article}
\newcommand{\Bbb}{\bf}
\def\frak#1{{#1}}

\def\ket#1{| #1\rangle}		

\def\be{\begin{equation}}
\def\en{\end{equation}}
\def\bea{\begin{eqnarray}}
\def\ena{\end{eqnarray}}
\def\bean{\begin{eqnarray*}}
\def\enan{\end{eqnarray*}}

\def\infq4p#1{{(#1;q^4,p)_\infty}}

\def\EXP#1{{\exp\biggl\{#1\biggr\}}}

\def\theta1{\vartheta_1}
\def\eh{\hat{e}_\lambda}
\def\fh{\hat{f}_\lambda}
\def\th{\hat{t}_\lambda}

\newcommand{\gsl}{{\frak {sl}}}

\def\slth{\widehat{\gsl_2}}
\def\aheslt{{\cal A}_{\hbar,\eta}(\slth)}
\def\aqpslt{{\cal A}_{q,p}(\slth)}
\def\uxpslt{U_{q,p}(\slth)}

\def\uqslth{U_q(\slth)}

\def\Ker{{\rm Ker}}
\def\Image{{\rm Im}}

\def\sh{{\rm sh}\,}

\def\tr{{\rm tr}\,}

\def\ket#1{|#1\rangle}

\def\Phil#1{\Phi^{(l)}_{#1}(\z)}
\def\Philz#1{\Phi^{(l)}_{#1}(z)}
\def\Psils#1{\Psi^{(l)*}_{#1}(\z)}
\def\Psilsz#1{\Psi^{(l)*}_{#1}(z)}

\def\F{{\cal F}}
\def\H{{\cal H}}
\def\Z{{\Bbb Z}}

\def\C{{\Bbb C}}
\def\R{{\Bbb R}}

\def\a{\alpha}
\def\b{\beta}

\def\z{\zeta}

\newtheorem{thm}{Theorem}[section]
\newtheorem{prop}[thm]{Proposition}
\newtheorem{lem}[thm]{Lemma}
\newtheorem{dfn}[thm]{Definition}
\newtheorem{conj}[thm]{Conjecture}
\newtheorem{cor}[thm]{Corollary}

\newcommand{\qed}{\hfill \fbox{}\medskip}
\newcommand{\proof}{\medskip\noindent{\it Proof.}\quad }

\newcommand{\refeq}[1]{(\ref{#1})}

\newcommand{\boszero}[2]{\, {\phi_0} ({#1}\, | \, {#2})\,}
\newcommand{\boszerop}[2]{\, {\phi'_0} ({#1}\, | \, {#2} )\,}
\newcommand{\bosone}[3]{\, {\phi_1}  \biggl(
{#1}\, \biggl| \, {#2} ; { #3} \biggr)\,}
\newcommand{\bostwo}[3]{\, {\phi_2}  \biggl(
{#1}\, \biggl| \, {#2} ; { #3} \biggr)\,}

\newcommand{\Ointz}[1]{\displaystyle \oint \frac{d {#1}}{2\pi i {#1}}\; }

\newcommand{\Oint}[1]{\displaystyle \oint \frac{d {#1}}{2\pi i}\; }
\newcommand{\OintC}[2]{\displaystyle \oint_{{#1}} \frac{d {#2}}{2\pi i}\; }

\newcommand{\zfactor}[2]{\, \frac{[{#1}]_{r-k}}{[{#2}]_{r-k}} \,}
\newcommand{\zfactorp}[2]{\, \frac{[{#1}]_r}{[{#2}]_r} \,}
\newcommand{\zfactorpf}[2]{\, \frac{[{#1}]_{k+2}}{[{#2}]_{k+2}} \,}

\newcommand{\zfactora}[2]{\, \frac{[{#1}]_{x}}{[{#2}]_{x}} \,}

\begin{document}
\font\csc=cmcsc10 scaled\magstep1

\begin{flushright}
September 1997 
\end{flushright}
\vspace{40pt}

\begin{center}
\begin{LARGE}
{ An Elliptic Algebra $U_{q,p}(\widehat{sl_2})$}\par
{and}\par
{the Fusion RSOS Model}
\vskip 3mm
\end{LARGE}

\vspace{30pt}

\begin{large}
Hitoshi KONNO
\end{large}

\vspace{6pt}

{\it Department of Mathematics, Faculty of Integrated Arts and Sciences,}\\
{\it  Hiroshima University, Higashi-Hiroshima 739,
 Japan\raisebox{2mm}{$\dagger$}}

\vspace{45pt}

{ABSTRACT}
\end{center}

\vspace{20pt}
We introduce an elliptic algebra
$U_{q,p}(\widehat{sl_2})$ with $p=q^{2r}\ (r\in \R_{>0})$ 
and present its  free boson representation 
at generic level $k$. 
We show that this algebra governs a structure of 
the space of states in the $k-$fusion RSOS model specified by a pair of 
positive integers $(r,k)$, or equivalently
 a $q-$deformation of the coset conformal field theory
$SU(2)_k\times SU(2)_{r-k-2}/SU(2)_{r-2}$.
Extending 
the work by Lukyanov and Pugai corresponding to the case $k=1$,
we gives 
a full set of screening operators for $k>1$.
The algebra $U_{q,p}(\widehat{sl_2})$ has two interesting degeneration limits, $p\to 0$
and $p\to 1$. The former limit yields the quantum affine algebra 
$U_{q}(\widehat{sl_2})$ whereas the latter
yields the algebra 
${\cal A}_{\hbar,\eta}(\widehat{sl_2})$,
the scaling limit of the elliptic algebra 
${\cal A}_{q,p}(\widehat{sl_2})$.
 Using this correspondence, we 
also obtain the highest component of 
two types of vertex operators which can be regarded as 
$q-$deformations of the primary fields in the coset conformal field 
theory.

\vspace{24pt}

\vfill
\hrule

\vskip 3mm
\begin{small}

\noindent\raisebox{2mm}{$\dagger$}E-mail: konno@mis.hiroshima-u.ac.jp
\end{small}

\setcounter{section}{0}
\setcounter{equation}{0}
\section{Introduction}
In studying exactly solvable models, especially in calculating correlation 
functions, the algebraic analysis method has proved to be extremely 
powerful\cite{JM}.
The method is based on the infinite 
dimensional quantum group symmetry possessed by a model\footnote{
Here we assume the model to have
 an infinite number of degrees of freedom.} and 
its representation theory. In particular, the intertwining operators
between the infinite dimensional representation spaces play an important role. 
There are two types of intertwiners called type I and type II. 
Remarkably, in solvable lattice models,
the type I intertwining operator can be identified with a 
certain composition of the Boltzmann weights so that 
 its product behaves as  
a local operator acting on the lattice. 
One can thus combine this into the Baxter's corner transfer matrix
(CTM) method\cite{Bax82}. 
As a consequence, the type I intertwiner allows us to identify the 
infinite dimensional irreducible representation with the space of states
of the model. Furthermore,
the  properties of the Boltzmann weights such as the 
Yang-Baxter equation, the inversion relation and the crossing 
symmetry yield some universal relations among the type I 
intertwiners and the CTM.  
Based on these relations, one can derive  $q-$difference equations
for correlation functions of local operators.

It is usually a difficult problem to solve such equations.
The advantage of the algebraic analysis 
is that it allows us to  
derive the correlation functions directly.
Correlation functions are formulated as 
traces of the product of type I intertwiners
over  irreducible representation space.   
Especially, in many cases, the infinite dimensional quantum group symmetries 
admit  free boson realization. 
This enables us to construct
the infinite dimensional representations as well as 
 the intertwining operators.
Then the calculation of the correlation functions i.e. the traces 
of the intertwiners 
is a straightforward task.
It is needless to say that the same spirit was 
already applied to the two dimensional conformal field 
theory(CFT)\cite{BPZ,DotFat}.

In \cite{JM,DFJMN,JMMN}, the $XXZ$ model, or equivalently, 
the six vertex model, 
in the anti-ferromagnetic regime 
was solved by applying the representation 
theory of the quantum affine algebra $\uqslth$. Following this work, 
its higher spin extension\cite{IIJMNT,Idzumi,BW,Konno94} and 
a higher rank extension\cite{Koyama} were discussed. The $XYZ$ model,
or equivalently, the eight vertex model in the principal regime, 
was also treated in this 
approach\cite{JMN,JKKMW}. 
There the elliptic algebra 
${\cal A}_{q,p}(\widehat{sl_2})$\cite{FIJKMY} was  
proposed  as the basic symmetry of the model. 
However its free field realizations still remain to be obtained.
It should also be remarked that the central extension of the Yangian double  
${\cal D}Y(\widehat{sl_2})$, the rational limit of $U_q(\widehat{sl_2})$, 
was properly defined\cite{Iohara,KLP96}. Its free field realization 
and application to physical problems were 
discussed\cite{Iohara,KLP96,Konno96}. 
 
In the recent work\cite{JKM}, we discussed  
two degeneration limits of the elliptic  algebra $\aqpslt$,
that are the limits $p\to 0$ and $p\to 1$. 
The first limit gives  
$\uqslth$, whereas the second limit yields a new algebra. 
The new algebra  turned out to be a relevant symmetry for 
the $XXZ$ model in the gap-less regime\cite{JKM} as well as for the 
sine-Gordon theory\cite{Lukyanov,Konno97}. 
This new algebra  was later reformulated  
by using  the Drinfeld currents and called 
$\aheslt$\cite{KLP}. 

On the other hand, it is known as vertex-face correspondence that there exists an interaction-round-a-face model corresponding to a vertex model.
The eight vertex model and the corresponding 
eight vertex solid-on-solid (SOS) model, 
or called  the Andrew-Baxter-Forrester (ABF) model, 
are  well-known examples\cite{ABF}.
The higher spin analogues of the ABF model were constructed by a fusion 
procedure\cite{DJKMO}.   
The $k-$fusion SOS model is obtained by fusing the Boltzmann weights 
of the ABF model $k$ times in both horizontal 
and vertical directions. 

We are interested in their 
restricted versions, i.e. the $k-$fusion restricted SOS (RSOS) models
\cite{ABF,DJKMO}.
The model is labeled by a pair of positive integers $(r,k)$.
At each lattice site $a$, one places a random variable (local height) $m_a$
taking values in the set $S=\{1,2,..,r-1\}$. One further imposes a restriction that for all adjacent sites $a$ and $b$, the local heights $m_a$ and $m_b$ 
satisfy the admissible conditions 
$$
m_a-m_b=-k,-k+2,..,k, \qquad
k+1\leq m_a+m_b \leq 2r-k-1.          
$$ 
The Boltzmann weight $W\left(\matrix{m_a&m_b\cr
                                     m_d&m_c\cr}\Bigl|u \right)$ 
of the model is attached to each configuration 
$(m_a,m_b,m_c,m_d)$ on the NW,NE,SE,SW corners of an elementary face
with $u$ being the spectral parameter.

The following two facts shown in \cite{DJKMO} for the regime III, 
which is treated through this paper, are fundamental.
\begin{enumerate}
\item The one point 
local height probability (LHP), i.e. 
a probability in which the center site has a given value of 
the local height, is given by the so-called  branching coefficient
appearing in the decomposition of the product of the two irreducible characters
of level $k$ and $r-k-2$ of the affine Lie algebra $A^{(1)}_1$ into 
level $r-2$ irreducible character of the diagonal $A^{(1)}_1$. 

\item The critical behavior is described 
by the conformal field theory with the Virasoro central charge 
$c=\frac{3k}{k+2}\Bigl(1-\frac{2(k+2)}{r(r-k)}\Bigr)$ 
and the primary fields of conformal dimensions 
$h_{J;n',n}$ \refeq{hw}.
\end{enumerate}

The corresponding CFT is known  as 
the coset minimal model
$SU(2)_k\times SU(2)_{r-k-2}/SU(2)_{r-2}$. The cases  $k=1,2$ were known 
before the fusion RSOS model\cite{GKO}, whereas
the  $k>2$ cases
were realized inspired by the model\cite{KMQ,BNY,Rava}.
The coset minimal model possesses the extended Virasoro algebra symmetry
generated by the Virasoro  generator and some extra generators.
The super Virasoro algebra is contained as the case $k=2$.
The irreducible characters of the extended Virasoro algebras 
are given by the same branching coefficients as the LHP in the above.

The first attempt to applying the algebraic analysis to 
the fusion RSOS model was carried out in \cite{JMOh}. 
There, in the regime III, the 
space of states was described based on 
the representation of $\uqslth$ and the creation and annihilation 
operators of the quasi-particles were obtained as a tensor product 
of certain type I and type II intertwiners in $\uqslth$. Later,
the type I vertex operator in the restricted ABF model in the regime III was 
recognized as a lattice operator and a $q-$difference equation for 
the correlation function was derived\cite{FJMMN}. 
   
Recently, Lukyanov and Pugai have succeeded to realize this type I
 vertex operator by using a free boson\cite{LukPug}.
 This enables us to construct 
a solution of the above $q-$difference equation exactly. However 
the most important 
contribution by this work is not this but a discovery of a symmetry 
of the restricted ABF model. 
That is the symmetry generated by 
the $q-$deformation of the Virasoro algebra\cite{AKOS}.  
In the same way as the Virasoro algebra, the $q-$Virasoro algebra admits a 
singular representation corresponding to the minimal series\cite{BPZ}. 
In such representation, 
screening operators play an essential role.
Constructing screening operators 
and $q-$deformation of the primary fields,
Lukyanov and Pugai discovered  
that their $q-$primary fields are nothing but the above type I vertex operator
in the restricted ABF model.

The purpose of this paper is to extend their result to the $k-$fusion RSOS
model in the regime III. 
Since the Virasoro algebra is realized as the case $k=1$ in the
coset CFT, we expect that a certain $q-$deformation of 
the extended Virasoro algebra corresponding to the coset
$SU(2)_k\times SU(2)_{r-k-2}/SU(2)_{r-2}$ $(k>1)$ exists and it provides  
a basic symmetry of the $k-$fusion RSOS models\cite{Konno97a}. 

Our strategy is based on the following 
observation. {\it 
The screening currents found by Lukyanov and Pugai satisfy an 
elliptic deformation of the quantum affine algebra $\uqslth$ at level one.  
This elliptic algebra has another degeneration limit to $\aheslt$ at 
level one\cite{JKM,AKOS}}. Picking up this algebraic nature, we
carry out the extension in the following two steps.
\begin{enumerate}
\item We  extend the elliptic algebra 
of the screening currents to generic level $k$. We call this algebra as 
$\uxpslt$\cite{Konno97a}. Realizing it by 
using free bosons, we show that the conformal limit of the   
generating currents for $\uxpslt$ 
coincide with those known in the above coset CFT. 
Hence these currents give a full 
set of screening currents for the $q-$deformed  coset theory
for arbitrary $k$. 
\item
The elliptic algebra $\uxpslt$ 
has two desired degeneration limits, $\uqslth$ and 
$\aheslt$. 
The Hopf algebra (like) structures
are known both in  $\uqslth$ and 
$\aheslt$. Using this knowledge, 
we obtain a free field 
realization of , at least, the highest component of the 
type I and type II vertex operators.
\end{enumerate}

The free field realization of 
the screening currents and the type I, II intertwiners
allows us to make a characterization of the 
$q-$deformation of the coset CFT $SU(2)_k\times SU(2)_{r-k-2}/SU(2)_{r-2}$. 

In order to 
identify our $q-$deformed coset theory 
with the $k-$fusion RSOS model, we investigate the following two things.
The first one is a partition function per site. 
We show that the correct partition function per site is obtained 
from the commutation of the two type I vertex operators. This allows us to
regard the type I vertex operator as a  proper lattice operator for
the $k-$ fusion RSOS model.
The second one is a 
structure of the Fock modules for the $q-$deformed  coset theory. 
The Fock modules are reducible due to the 
existence of singular vectors which can be constructed 
by the screening operators on the modules.
We  consider a resolution of the modules  
and obtain a space which can be regarded as the  irreducible highest 
weight representation of the conjectural $q-$deformation of the extended 
Virasoro algebra. 
The character 
of the space coincides with the desired branching 
coefficient. Hence one  can identify the space with the  
space of states for the $k-$fusion RSOS model. 
   
This paper is organized as follows. In the next section, we briefly review the 
free field representation of the 
coset minimal model $SU(2)_k\times SU(2)_{r-k-2}/SU(2)_{r-2}$. 
The formulae summarized in this section 
should be compared with those of the $q-$deformed ones obtained in Sec.4 and 5.
In Sec.3, we introduce the elliptic algebra $\uxpslt$ and discuss its 
properties. In Sec.4, we present a free field representation of 
$\uxpslt$. As a corollary, the free field representation
of the algebra  $\aheslt$ for arbitrary level is obtained. 
We also derive  the type I and type II vertex operators
and their commutation relations for the highest components.     
We argue that  the correct partition function per site is obtained from 
these relations.
In Sec.5, based on these results, we propose a $q-$deformation of the 
coset minimal model. We show that the algebra 
$\uxpslt$ provides a full set of screening
operators which are sufficient for making 
a resolution of the Fock 
modules. Then a characterization of the 
irreducible highest weight modules of the
conjectural $q-$deformed extended Virasoro algebras is obtained. 
The evaluation of the character of the space supports the 
identification of the space with the space of states for  
the $k-$fusion RSOS model.
The final section is devoted to discussions on the results
and some future problems.

\setcounter{section}{1}
\setcounter{equation}{0}
\section{Coset Conformal Field Theory}
      In this section we  briefly review the  
coset minimal
model $SU(2)_k\times SU(2)_l/SU(2)_{k+l}$ 
$k,l\in\Z_{>0}$\cite{GKO,KMQ,BNY,Rava}.    

The symmetry of the theory is an extended Virasoro algebra 
 generated by the Energy-Momentum(EM) tensor $T(z)$ and the extra generators 
$A_k^j(z)\ (j=1,2,..)$. 
The number of extra generators  is depend on the value $k$.
 For example, the case $k=1$ and $l\in\Z_{>0}$, the theory 
is the Virasoro minimal model and there are no extra generators\cite{BPZ,GKO}. 
 For $k=2$
with $l\in\Z_{>0}$, the theory is the 
$N=1$ super Virasoro minimal model\cite{GKO,BKT,FQS}. 
There are one  extra generator of conformal dimension 
$3/2$, which is nothing but the super generator $A_2(z)=G(z)$.  
The $k=4$ case with $l\in\Z_{>0}$ is known as the $S_3$ 
symmetric model\cite{ZamFat}.
There are two  extra generators $A_4^1(z)$ and $A_4^2(z)$
 of conformal dimension $4/3$. 

The extended Virasoro algebra is defined by the following operator product expansions (OPE)\cite{Rava}.
\bea
&&T(z)T(w)=\frac{c/2}{(z-w)^4}+\frac{2T(z)}{(z-w)^2}+
\frac{\partial T(z)}{z-w}+O(1)\\
&&T(z)A^j_k(w)=\frac{\sigma A^j_k(w)}{(z-w)^2}+
\frac{\partial A^j_k(w)}{z-w}+O(1)\\
&&A^i_k(z)A^j_k(w)=\alpha^{ij}(z-w)^{-2\sigma}\Bigl\{
\frac{c}{\sigma}+2T(w)(z-w)^2+O((z-w)^3)\Bigr\}\\
&&\qquad +\beta^{ijm}(z-w)^{-\sigma}\Bigl\{
A^m_k(w)+\frac{1}{2}(z-w)\partial A^m_k(w)+O((z-w)^2)\Bigr\},  
\ena
where $\sigma=\frac{k+4}{k+2}$ is the conformal dimension of
 $A_k^j(z)$, and $\a^{ij}, \b^{ijm}$ are the structure constants.
The central charge $c$ of the Virasoro algebra is given by
\bea
c=\frac{3k}{k+2}\Bigl(1-\frac{2(k+2)}{(l+2)(l+k+2)}\Bigr).
\label{centralc}
\ena

Let 
$r=l+k+2$, $1\leq n\leq r-k-1,\ 1\leq n'\leq r-1$
and $J=|n'-n ({\rm mod} 2k)|,\  0\leq J\leq k$.
Define the Virasoro generators $L_m$ by $T(z)=\sum_{m\in \Z}
L_m z^{-m-2}$.
The Virasoro highest weight state $\ket{ J;n',n}$ 
is the state  satisfying
\bea
&&L_0\ket{ J;n',n}=h_{J;n',n}\ket{ J;n',n},\\
&&L_n\ket{ J;n',n}=0,  \qquad {\rm for }\qquad n\geq1.
\label{vhws}
\ena
The highest weights are given by
\bea
&&h_{J;n'n}={J(k-J)\over 2k(k+2)}+{(nr-n'(r-k))^2-k^2\over 4kr(r-k)}.
\label{hw}
\ena

One can realize the theory in terms of three boson 
fields\cite{KMQ,BNY,GMM}.
Let us introduce the 
 three independent 
free boson fields $\phi_0(z), \phi_1(z)$ and $ \phi_2(z)$
satisfying the OPE $<\phi_0(z)\phi_0(w)>=<\phi_1(z)\phi_1(w)>=
\log(z-w)=-<\phi_2(z)\phi_2(w)>$. 
Then the  EM tensor of the coset theory is realized as 
\bea
&&T(z)=T_{Z_k}(z)+
{1\over 2}(\partial\phi_0(z))^2+
\sqrt{2}\alpha_0\partial^2 \phi_0(z) ,
\label{emcoset}\\
&&T_{Z_k}(z)=
{1\over 2}(\partial \phi_1(z))^2-{1\over 2}\sqrt{\frac{2}{k+2}}
\partial^2 \phi_1(z)
-{1\over 2}(\partial \phi_2(z))^2.
\label{empara}
\ena
Here  $2\alpha_0=\sqrt{k\over r(r-k)}, $ and 
$T_{Z_k}(z)$ being  the EM tensor of the $Z_k$ parafermion theory\cite{ZamFat}.
The realization of the extra generators can be found 
in \cite{KMQ,BNY,Rava}.

The expression \refeq{emcoset} indicates that
the coset theory is realized as a composition of the $\Z_k$ 
parafermion theory and the  boson theory $\phi_0(z)$. 
Indeed,
the primary field of the coset theory is realized as
\bea
&&\Psi_{J;n',n}(z)=
\Psi_{J,M;n',n}(z)|_{J=M ({\rm mod}\ 2k)},\\
&&\Psi_{J,M;n',n}(z)=\phi_{{J},M}(z):\exp
\sqrt2\alpha_{n',n}\phi_0(z):,\label{primaryf}\\
&&\phi_{J,M}(z)
=:\exp\Bigl\{\frac{J}{\sqrt{{2}{k+2}}}\phi_1(z)+ \frac{M}{\sqrt{{2}{k}}}
\phi_2(z)\Bigr\}:,\label{zkprimaryf}
\ena
where  $M=-J,-J+2,..,J$ (mod $2k$), 
$\alpha_{n',n}={1-n'\over 2}\alpha_{-}+{1-n\over 2}\alpha_{+},\ 
\alpha_+=\sqrt{\frac{r}{k(r-k)}},\ \alpha_-=-\sqrt{\frac{r-k}{kr}}$.
Note $2\a_0=\a_++\a_-$.
Using $\Psi_{J;n',n}(z)$, the highest weight state $\ket{J;n',n}$ is obtained  as
\bea
&&\ket{J;n',n}=\lim_{z\to 0}\Psi_{J;n',n}(z)\ket{0},
\label{hwstate}
\ena
where $\ket{0}$ denotes the $SL(2,\C)$ invariant vacuum state.

The highest weight representation space is then given by  
the Fock module  ${\cal F}_{J;n',n}$
constructed  by the 
action of the creation operators of the fields $\phi_j(z), \ j=0,1,2$ 
on $\ket{ J;n',n}$. These modules are reducible due to the existence of 
singular vectors. 
The singular vectors can be constructed by using the screening operators
on some highest weight states.
In the minimal coset theory, the screening currents, 
which contour integrals yield the  screening operators,
are given by
\bea
&&S_{+}(z)=\Psi(z) :\exp\Bigl\{\sqrt 2\alpha_{+}\phi_0(z)\Bigr\}:, 
\label{scp}\\ 
&&S_{-}(z)=\Psi^{\dagger}(z) :\exp\Bigl\{\sqrt 2\alpha_{-}\phi_0(z)\Bigr\}:,
\label{scm}\\
&&S(z)=-:\Bigl(\sqrt{\frac{k+2}{2}}\partial\phi_1(z)+
\sqrt{\frac{k}{2}}\partial\phi_2(z)\Bigr)\exp
\Bigl\{-\sqrt{\frac{2}{k+2}}\phi_1(z)\Bigr\}:,
\\
&&\eta(z)=:\exp\Bigl\{\sqrt{\frac{k+2}{2}}\phi_1(z)+
\sqrt{\frac{k}{2}}\phi_2(z)\Bigr\}:,
\ena 
where
$\Psi$ and $\Psi^{\dagger}$ are the parafermion currents given 
by
\bea
&&\Psi(z)=:\Bigl(\sqrt{\frac{k+2}{k}}\partial\phi_1(z)+\partial\phi_2(z)\Bigr)\exp\Bigl\{\sqrt{2\over k}\phi_2(z)\Bigr\}:,\\
&&\Psi^{\dagger}(z)=:\Bigl(\sqrt{\frac{k+2}{k}}\partial\phi_1(z)-\partial\phi_2(z)\Bigr)
\exp\Bigl\{-\sqrt{2\over k}\phi_2(z)\Bigr\}:.
\ena
These currents are characterized by the properties 
(i) the conformal dimension is one, (ii) 
the contour integrals of them commute 
with the extended Virasoro algebra. 

These screening currents 
defines the nilpotent operators called the BRST charges.
One can use these charges to make a resolution of the Fock modules\cite{
Felder,BerFel}.
One should note that the screening currents 
$S(z)$ and $\eta(z)$ act only on the $\Z_k$ parafermion theory and that 
the screening operators obtained from $S_\pm(z)$ commute with those from
$S(z)$ and $\eta(z)$.
Therefore, in order to make a resolution of 
the Fock modules of the coset theory, 
one may take the following two steps. First make a resolution of  the 
$\Z_k$ parafermion theory and then consider the coset 
theory\cite{Konno92}.\footnote{
The boson fields $\Phi(z), \varphi(z), \chi(z)$ 
used in \cite{Konno92} are related to 
$\phi_0(z), \phi_1(z), \phi_2(z)$ as follows.
\bean
&&\Phi(z)=\phi_0(z),\quad\varphi(z)=-\sqrt{\frac{k}{2}}\phi_1(z)-
\sqrt{\frac{k+2}{2}}\phi_2(z),\quad
\chi(z)=-\sqrt{\frac{k+2}{2}}\phi_1(z)-
\sqrt{\frac{k}{2}}\phi_2(z).
\enan} 
In section 5, we discuss a $q-$analogue of this  resolution.

\setcounter{section}{2}
\setcounter{equation}{0}
\section{The Elliptic Algebra $U_{q,p}(\widehat{sl_2})$}
In this section we define a new elliptic
 algebra $U_{q,p}(\widehat{sl_2})$\cite{Konno97a} and 
discuss its relation to known algebras $\uqslth$ and 
$\aheslt$. 

\subsection{Definition}
Let $r\in \R_{>0}$ and  $q\in \C,\ |q|<1$.
We set $p=q^{2r}$ and $p^*=p q^{-2c}$.

\begin{dfn}
The associative algebra $U_{q,p}(\widehat{sl_2})$ is generated by 
the operator valued functions (currents) 
$k(z), \ E(z), \ F(z)$ of complex variable $z$ and the 
central element $c$ 
with the following relations\footnote{We are indebted to Jimbo for 
introducing the generator $k(z)$ and giving its relation to 
$H^{\pm}(z)$ in \refeq{hkk}. 
}.
\bea
&&k(z)k(w)=\Bigl(\frac{z}{w}\Bigr)^{\frac{c}{2r(r-c)}}
\frac{\xi(w/z;p,q)}{
\xi(w/z;p^*,q)}
\frac{\xi(z/w;p^*,q)}{
\xi(z/w;p,q)}k(w)k(z),\label{defa}\\
&&k(z)E(w)=
\Bigl(\frac{z}{w}\Bigr)^{\frac{1}{r-c}}
\frac{\Theta_{p^*}(q^{-1}p^{*\frac{1}{2}}w/z)}{
\Theta_{p^*}(q p^{*\frac{1}{2}}w/z)}
E(w)k(z),\label{defb}\\
&&k(z)F(w)=
\Bigl(\frac{z}{w}\Bigr)^{-\frac{1}{r}}
\frac{\Theta_{p}(q p^{\frac{1}{2}}w/z)}{
\Theta_{p}(q^{-1} p^{\frac{1}{2}}w/z)}
F(w)k(z),\\
&&E(z)E(w)=
q^2\Bigl(\frac{z}{w}\Bigr)^{\frac{2}{r-c}}
\frac{\Theta_{p^*}(q^{-2}w/z)}{
\Theta_{p^*}(q^{2}w/z)}
E(w)E(z),\\
&&F(z)F(w)=
q^{-2}\Bigl(\frac{z}{w}\Bigr)^{-\frac{2}{r}}
\frac{\Theta_{p}(q^{2}w/z)}{
\Theta_{p}(q^{-2}w/z)}
F(w)F(z),\label{deff}\\
&&[E(z),F(w)]=\frac{1}{q-q^{-1}}
\Bigl\{\delta(q^{-c}z/w)H^+(q^{-\frac{c}{2}}z)
-\delta(q^{c}z/w)H^-(q^{-\frac{c}{2}}w)\Bigr\},
\label{defg}
\ena
where 
\bea
&& 
H^{\pm}(z)=\kappa(z) k(q^{\pm(r-\frac{c}{2})+1}z)
k(q^{\pm(r-\frac{c}{2})-1}z),\label{hkk}\\
&& \kappa(z)=\Bigl(q^{\pm(r-\frac{c}{2})+1}z\Bigr)^{-\frac{c}{2r(r-c)}}
\frac{\xi(x;p^*,q)}{\xi(x;p,q)}\Bigr|_{x=q^{-2}},\\
&&\xi(z;p,q)=\frac{(q^2 z;p,q^4)_{\infty}(p q^2 z;p,q^4)_{\infty}}
{(q^4 z;p,q^4)_{\infty}(p z;p,q^4)_{\infty}}
\ena
and 
$[A,B]=AB-BA$, $\delta(z)=\sum_{n\in\Z}z^{n}$,
$\Theta_s(z)=(z;s)_\infty(s/z;s)_\infty(s;s)_\infty ,\ $ 
$(z;s)_\infty=\prod_{n=0}^\infty(1-z s^n)$
\end{dfn}

\vspace{3mm}
\noindent
{\it Remark 3.1.}
Let us set $p=e^{-2\pi i/\tau},\ p^*=e^{-2\pi i/\tau^*}$, $z=q^{2u}$ and  
denote $E(z), \ F(z)$ and $k(z)$ by the same letters $E(u),\ F(u)$ and
$k(u)$.
Then the relations in (\ref{defb})-\refeq{deff} are rewritten in
more compact form.
\bea
&&k(u)E(v)=
\frac{\theta1(\frac{u-v+\frac{1}{2}}{r-c}-\frac{1}{2}|\tau^*)}
{\theta1(\frac{u-v-\frac{1}{2}}{r-c}-\frac{1}{2}|\tau^*)}
E(v)k(u),\\
&&k(u)F(v)=
\frac{\theta1(\frac{u-v-\frac{1}{2}}{r}-\frac{1}{2}|\tau)}
{\theta1(\frac{u-v+\frac{1}{2}}{r}-\frac{1}{2}|\tau)}F(v)k(u),\\
&&E(u)E(v)=
\frac{\theta1(\frac{u-v+1}{r-c}|\tau^*)}
{\theta1(\frac{u-v-1}{r-c}|\tau^*)}E(v)E(u),\\
&&F(u)F(v)=
\frac{\theta1(\frac{u-v-1}{r}|\tau)}
{\theta1(\frac{u-v+1}{r}|\tau)}F(v)F(u), \\
&&H^{+}(u)H^{-}(v)=\frac{\theta1(\frac{u-v-(1+\frac{c}{2})}{r}|\tau)}{
\theta1(\frac{u-v+(1-\frac{c}{2})}{r}|\tau)}
\frac{\theta1(\frac{u-v+(1+\frac{c}{2})}{r-c}|\tau^*)}{
\theta1(\frac{u-v-(1-\frac{c}{2})}{r-c}|\tau^*)}H^{-}(v)H^{+}(u),\\
&&H^{\pm}(u)H^{\pm}(v)=\frac{\theta1(\frac{u-v-1}{r}|\tau)}{
\theta1(\frac{u-v+1}{r}|\tau)}
\frac{\theta1(\frac{u-v+1}{r-c}|\tau^*)}{
\theta1(\frac{u-v-1}{r-c}|\tau^*)}H^{\pm}(v)H^{\pm}(u),
\ena
where  
$\theta1(u|\tau)$ is the Jacobi elliptic theta function  
\bean
\theta1(u|\tau)=i\sum_{n\in \Z}(-)^ne^{\pi i (n-1/2)^2\tau}
e^{2\pi i (n-1/2) u}.
\enan
These expressions suggest that the algebra $\uxpslt$ is related to some 
elliptic curves in the similar way to the theory of Enriquez and 
Felder\cite{EnriFeld}.


\subsection{Degeneration limit}
There are two interesting degeneration limits: $p\to 0$ and $p\to 1$.

The limit  $p\to 0$ is taken by letting  $ r\to \infty $. Then 
the relations (\ref{defa})-\refeq{defg} are reduced to 
those of the Drinfeld currents in the 
quantum affine algebra $U_q(\widehat{sl_2})$ ( see for example,
\cite{FreJin}).
In the later section (Sec.5), we will make the  identification that  
$U_{q,p}(\widehat{sl_2})$ is the algebra of the screening currents 
in the $q-$deformation of the coset theory 
$SU(2)_k\times SU(2)_{r-k-2}/SU(2)_{r-2}$. 
The limit $r\to \infty$ to $\uqslth$ is then 
consistent with the well-known fact  in CFT and 
the perturbation of the coset CFT\cite{ABL}, i.e.  
$$
\lim_{r\to \infty} SU(2)_k\times SU(2)_{r-k-2}/SU(2)_{r-2}\ {\rm theory}
\to SU(2)_k\ {\rm  Wess-Zumino-Witten\  model}.
$$   

The second limit $p\to 1$ is taken by setting 
$q=e^{\frac{\hbar \varepsilon}{2}}$ and $z=e^{-i\alpha\varepsilon},
w=e^{-i\beta\varepsilon}$ and letting $\varepsilon \to 0$. In this limit,
the relations in Definition 3.1 tend to  those of the currents in 
the algebra $\aheslt$ \cite{KLP} (see Appendix), i.e.
the degeneration (or scaling) limit of the  
elliptic algebra ${\cal A}_{q,p}(\widehat{sl_2})$\cite{FIJKMY}, 
 under the identification
$1/\eta=\hbar r, 1/\eta'=\hbar(r-k)$ and the interchanges 
$H^+\leftrightarrow H^-,\ E\leftrightarrow F$. 
This limit is also consistent with the known fact 
that the scaling limit
of the RSOS model gives the restricted sine-Gordon model
and the latter model is obtained by an integrable perturbation of the
coset minimal model(see, for example,\cite{ABL}).

Our algebra $\uxpslt$ hence has the same degeneration limits as the 
elliptic algebra ${\cal A}_{q,p}(\widehat{sl_2})$\cite{JKM}. 
The direct relationship between $\uxpslt$ and ${\cal A}_{q,p}(\widehat{sl_2})$
is unknown at this stage\footnote{In the trigonometric limit, 
the corresponding problem has been discussed by Hou et al.\cite{HouYang}.}.
However, see \cite{JKM} Sec.3 where one can find a discussion which
suggests the deep relation between the $q-$Virasoro
algebra and the elliptic algebra ${\cal A}_{q,p}(\widehat{sl_2})$
at level one.

\setcounter{section}{3}
\setcounter{equation}{0}
\section{Free Field Realization of $U_{q,p}(\widehat{sl_2})$ at Level k}
We here consider a realization of the algebra $U_{q,p}(\widehat{sl_2})$
at arbitrary level $k\not=0,-2$ in terms of three bosonic fields\footnote{
In preparing this paper, we have noticed that Shiraishi has obtained 
another free field realization for the similar algebra.}.

\subsection{Bosonization of $\uxpslt$}
Let $a_{j,m}\ ( m\in\Z_{\not=0}\ j=0,1,2)$ be bosons satisfying 
the commutation relations
\bea
&&[a_{0,m},a_{0,n}]=\frac{[2m][km]}{m}\frac{[rm]}{[(r-k)m]}\delta_{m+n,0},
\label{ccra}\\
&&[a_{1,m},a_{1,n}]=\frac{[2m][(k+2)m]}{m}\delta_{m+n,0},\label{ccrb}\\
&&[a_{2,m},a_{2,n}]=-\frac{[2m][km]}{m}\delta_{m+n,0},\label{ccrc}
\ena
where $[m]=\frac{q^m-q^{-m}}{q-q^{-1}}$.
We also define the primed boson $a'_{0,m}$
and the zero-mode operators $Q_j$ and $P_j\ (j=0,1,2)$ satisfying
\bea
&&a'_{0,m}=\frac{[(r-k)m]}{[rm]}\ a_{0,m},\\
&&[P_0,Q_0]=-i,\qquad [P_1,Q_1]=2(k+2),\qquad [P_2,Q_2]=-2k.
\ena

In order to make the expression of the currents simple, we introduce 
the following 
boson fields $\phi_j(A;B,C|z;D)$  $(j=0,1,2,\quad  A,B,C,D\in \R)$.
\bea
&&\phi_0(A;B,C|z;D)=\frac{A}{BC}\sqrt{\frac{2kr}{r-k}}
(i Q_0+P_0\log z)+\sum_{{m\in\Z \atop \not=0}}\frac{[Am]}{[Bm][Cm]}a_{0,m}z^{-m}
q^{D|m|},\label{bosonf0}\\
&&\phi'_0(A;B,C|z;D)=\frac{A}{BC}\sqrt{\frac{2k(r-k)}{r}}
(i Q_0+P_0\log z)+\sum_{{m\in\Z\atop \not=0}}\frac{[Am]}{[Bm][Cm]}a'_{0,m}z^{-m}
q^{D|m|},\nonumber\\
\\
&&\phi_j(A;B,C|z;D)=-\frac{A}{BC}
(Q_j+P_j\log z)+
\sum_{{m\in\Z\atop \not=0}}\frac{[Am]}{[Bm][Cm]}a_{j,m}z^{-m}
q^{D|m|}\quad (j=1,2),\nonumber\\
\\
&&\phi_j^{(\pm)}(A;B|z;C)=\frac{P_j}{2}\log q+
(q-q^{-1})\sum_{m\in\Z_{>0}}\frac{[Am]}{[Bm]}a_{j,\pm m}z^{\mp m}
q^{C m}\qquad (j=1,2).\nonumber\\
\label{bosonfjpm}\ena
We often use the abridgment
\bea
&&\phi_j(C|z;D)=\phi_j(A;A,C|z;D),\qquad \phi_j(C|z)=\phi_j(C|z;0). 
\ena
We denote by :: the usual normal ordered product. For example,
\bea
&&:\EXP{-\phi_0(k|z)}:=\EXP{\sum_{m\in \Z_{>0}}\frac{a_{0,-m}}{[km]}z^m}
\EXP{-\sum_{m\in \Z_{>0}}\frac{a_{0,m}}{[km]}z^{-m}}\nonumber\\ 
&&\qquad\qquad\qquad\qquad\qquad\times e^{-i\sqrt{\frac{2r}{k(r-k)}}Q_0}\ 
z^{-\sqrt{\frac{2r}{k(r-k)}}P_0}.
\ena

Then we have 
\begin{thm}
The algebra $U_{q,p}(\widehat{sl}_2)$ has the following free field realization at $c=k$.
\bea
&&k(z)=:\EXP{-\phi_0(1;2,r|z)}:,\label{currk}\\
&&E(z)= 
\Psi(z)\ :\EXP{-\phi_0(k|z)}:,\label{curre}\\
&&F(z)=
\Psi^{\dagger}(z)\ :\EXP{\phi'_0(k|z)}:\label{currf},
\ena
where $\Psi(z)$ and $\Psi^{\dagger}(z)$ are the $q-$deformed 
$\Z_k$ parafermion currents given by $\Psi(z)=\Psi^-(z),\ 
\Psi^{\dagger}(z)=\Psi^+(z) $\cite{Matsuo}\footnote{
The relation of our notations to those of Matsuo's is as follows.
$
\beta_m=a_{1,m}, \ \bar{\a}_m=a_{2,m},\ \b_0=P_1,\ \bar\a_0=P_2,\ 
2(k+2)\b=Q_1,\ 2k\bar\a=Q_2.
$
},
\bea
&&\Psi^{\pm}(z)=\mp\frac{1}{(q-q^{-1})}
:\EXP{\pm\phi_2(k|z;\pm\frac{k}{2})}\nonumber\\
&&\qquad\qquad\qquad\times\Bigl(
\EXP{-\phi^{(+)}_2\Bigl(1;2|z;\mp\frac{k+2}{2}\Bigr)\pm
\phi^{(+)}_1\Bigl(1;2|z;\mp\frac{k}{2}\Bigr)}\nonumber\\
&&\qquad\qquad\qquad\qquad
 -\EXP{\phi^{(-)}_2\Bigl(1;2|z;\mp\frac{k+2}{2}\Bigr)\mp
\phi^{(-)}_1\Bigl(1;2|z;\mp\frac{k}{2}\Bigr)}\Bigr):\label{pfcurrents}.
\ena
\end{thm}

\proof
Straightforward calculation.
\qed

\vspace{3mm}
\noindent
{\it Remark 4.1} In the CFT limit $q\to 1$ with $z$ fixed, the currents 
$E(z)$ and  $F(z)$ coincide with the screening currents $S_+(z)$\refeq{scp}, 
$S_-(z)$\refeq{scm} in the coset CFT. 

\vspace{3mm}
\noindent
{\it Remark 4.2} In the limit $r\to \infty$, the expressions \refeq{currk}
with \refeq{hkk} and \refeq{curre}-
\refeq{currf} tend
to the Matsuo's bosonization of $U_q(\widehat{sl_2})$\cite{Matsuo}.
Namely, $H^+(z)\to K_+(z), \ H^-(z)\to K_-(z), E(z)\to X^+(z)$ and 
$F(z)\to X^-(z)$.  

\subsection{The ${\cal A}_{\hbar,\eta}(\widehat{sl_2})$ limit}
Let us consider 
the ${\cal A}_{\hbar,\eta}(\widehat{sl_2})$ limit.
The algebra $\aheslt$\cite{KLP} is defined in Appendix.
The limit  is 
taken by the following procedure\cite{JKM}. 
Setting $q=e^{\frac{\hbar \varepsilon}{2}}$, $z=e^{-i\alpha\varepsilon}$,
$r=\xi+k$, $m\varepsilon =t\in\R$ and letting $\varepsilon \to 0$. Then we have
from \refeq{ccra}-\refeq{ccrc}
\bea
&&[a_{0}(t),a_{0}(t')]=\frac{1}{\hbar^2}
\frac{\sinh\hbar t\sinh\frac{\hbar k t}{2}
\sinh\frac{\hbar(\xi+k)t}{2}}{t\sinh\frac{\hbar \xi t}{2}}\delta(t+t'),\\
&&[a_{1}(t),a_{1}(t')]=\frac{1}{\hbar^2}
\frac{\sinh\hbar t\sinh\frac{\hbar( k+2) t}{2}}{t}\delta(t+t'),\\
&&[a_{2}(t),a_{2}(t')]=-\frac{1}{\hbar^2}
\frac{\sinh\hbar t\sinh\frac{\hbar k t}{2}}{t}\delta(t+t'),
\ena
and
\bea
&&a'_{0}(t)=
\frac{\sinh\frac{\hbar \xi t}{2}}{\sinh\frac{\hbar (\xi+k) t}{2}}\ a_{0}(t).
\ena
In this limit, all the zero-mode operators $Q_j$ and $P_j\ (j=0,1,2)$ 
are dropped.

The boson fields $\phi_j(A;B,C|z;D),\ j=0,1,2$ are now reduced to 
\bea
&&\tilde\phi_j(A;B,C|\alpha;D)=\hbar\int_{-\infty}^{\infty}dt
\frac{\sinh\frac{\hbar At}{2}}{\sinh\frac{\hbar Bt}{2}\sinh\frac{\hbar Ct}{2}}
a_{j}(t)e^{i\alpha t+\frac{\hbar D}{2}|t|}\ (j=0,1,2),\nonumber\\
\\
&&\tilde\phi'_0(A;B,C|\alpha;D)=\hbar\int_{-\infty}^{\infty}dt
\frac{\sinh\frac{\hbar At}{2}}{\sinh\frac{\hbar Bt}{2}\sinh\frac{\hbar Ct}{2}}
a'_{0}(t)e^{i\alpha t+\frac{\hbar D}{2}|t|}
,\\
&&\tilde\phi_j^{(\pm)}(A;B|\alpha;C)=2\hbar\int_{0}^{\infty}dt
\frac{\sinh\frac{\hbar At}{2}}{\sinh\frac{\hbar Bt}{2}}
a_{j}(\pm t)e^{\pm i\alpha t+\frac{\hbar C}{2}|t|}\qquad (j=1,2).
\ena

Under these notations, we have

\begin{thm}
 The level $k$ currents in $\aheslt$ with 
$ 1/\eta=\hbar\xi,\ 1/\eta'=\hbar(\xi+k)$ is realized as follows.
\bea
&&k(\a)=:\EXP{-\tilde\phi_0(1;2,\xi+k\Bigr|\alpha)}:,\\
&&E(\alpha)=\Psi(\a)\ :\EXP{-\tilde\phi_0(k|\alpha)}:,\\
&&F(\alpha)= \Psi^{\dagger}(\a)\ :\EXP{\tilde\phi'_0(k|\alpha)}:,
\ena
where $\Psi(\alpha)=\Psi^-(\a),\ \Psi^{\dagger}(\a)=\Psi^+(\a) $ are given by
\bea
&&\Psi^{\pm}(\a)=\mp\frac{1}{\hbar}
:\EXP{\pm\tilde\phi_2\Bigl(k\Bigr|\a;\pm\frac{k}{2}\Bigr)}\nonumber\\
&&\times\Bigl(
\EXP{-\tilde\phi^{(+)}_2\Bigl(1;2\Bigr|\a;\mp\frac{k+2}{2}\Bigr)\pm
\tilde\phi^{(+)}_1\Bigl(1;2\Bigr|\a;\mp\frac{k}{2}\Bigr)}\nonumber\\
&&\qquad -\EXP{\tilde\phi^{(-)}_2\Bigl(1;2\Bigr|\a;\mp\frac{k+2}{2}\Bigr)\mp
\tilde\phi^{(-)}_1\Bigl(1;2\Bigr|\a;\mp\frac{k}{2}\Bigr)}\Bigr):.
\ena
\end{thm}

In the algebra $\aheslt$, the Hopf algebra like structure and 
the level zero representation are known\cite{KLP}. 
We summarize them in Appendix. 
From these knowledge, one can construct the 
intertwining operators between  the level $k$
infinite dimensional representations. 
Let us define the following four vertex operators.
\bea
&&\Phil{l}=\phi_{l,l}(\z):\EXP{-\tilde\phi_0'(l;2,k|\z)}:,\label{dtypeI}\\
&&\Psils{l}=\phi_{l,-l}(\z):
\EXP{\tilde\phi_0(l;2,k|\z)}:\label{dtypeII},\\
&&\tilde\Phil{l}=\phi_{k-l,-(k-l)}(\z):\EXP{-\tilde\phi_0'(l;2,k|\z)}:,
\label{dtypeIt}\\
&&\tilde\Psils{l}=\phi_{k-l,k-l}(\z):
\EXP{\tilde\phi_0(l;2,k|\z)}:\label{dtypeIIt}\qquad (l=0,1,2,..,k),
\ena
where $\phi_{l,\pm l}(\z)$ are the analogues of the $\Z_k$-parafermion
primary fields \refeq{zkprimaryf} given by: 
\bea
&&\phi_{l,\pm l}(\z)=
:\EXP{-\tilde\phi_2\Bigl(\pm l;2,k|\z;\pm\frac{k}{2}\Bigr)-
\tilde\phi_1\Bigl(l;2,k+2|\z;\pm\frac{k+2}{2}\Bigr)}:.
\ena
Then we have

\begin{prop}
The vertex operators $\Phil{l}$ and $\Psils{l}$ satisfy the intertwining 
relations of the type I \refeq{inthi}-\refeq{intfi} and the type II 
\refeq{inthii}-\refeq{inteii}, respectively, whereas
the vertex operators $\tilde\Phil{l}$ and $\tilde\Psils{l}$ 
satisfy the twisted intertwining 
relations of the type I \refeq{inthi}, \refeq{tintei}-\refeq{tintfi} 
and the type II \refeq{inthii},\refeq{tintfii}-\refeq{tinteii}, respectively. 
\end{prop}
\proof
Straightforward.
\qed

Using the relations \refeq{intfi}-\refeq{inteii},
one can obtain the other (lower) components of the intertwiners 
$\Phil{m},\ \Psils{m},\ \tilde\Phil{m},\ \tilde\Psils{m}\ (m=0,1,2,..,l-1)$.
We omit them here.

\subsection{The type I and type II vertex operators}

Although we do not have any results on the Hopf algebra structure of
$\uxpslt$, the $(q,p)$-analogue of the 
operators \refeq{dtypeI}-\refeq{dtypeIIt}
 can be obtained by the following requirements.

\begin{enumerate}
\item The procedure taking the limit to $\aheslt$
from $\uxpslt$ makes the vertex operators 
reduce to those in \refeq{dtypeI}--\refeq{dtypeIIt}. 
\item The zero-modes factors are determined by 
 requiring that the CFT limit
( $q\to 1$ with $z$ fixed ) of the vertex operators should be expressed 
as the exponential of the boson fields.  
\end{enumerate}

Using the notations in \refeq{bosonf0}-\refeq{bosonfjpm}, we find that the desired vertex operators 
are given as follows.
\bea
&&\Philz{l}=\phi_{l,l}(z):
\EXP{-\phi_0'(l;2,k|z)}:,\label{typeivo}\\
&&\Psilsz{l}=\phi_{l,-l}(z):
\EXP{\phi_0(l;2,k|z)}:,\label{tupeiivo}\\
&&\tilde\Philz{l}=\phi_{k-l,-(k-l)}(z):
\EXP{-\phi_0'(l;2,k|z)}:,\label{ttypeivo}\\
&&\tilde\Psilsz{l}=\phi_{k-l,k-l}(z):
\EXP{\phi_0(l;2,k|z)}:\qquad(l=0,1,2,..,k),\label{ttypeiivo}
\ena
where
\bea
&&\phi_{l,\pm l}(z)=:\EXP{-\phi_2\Bigl(\pm l;2,k|z;\pm\frac{k}{2}\Bigr)-
\phi_1\Bigl(l;2,k+2|z;\pm\frac{k+2}{2}\Bigr)}:.
\ena

We also have some conjectural expression for the lower components of these 
vertex operators. We will discuss them and their commutation relations 
in a separate paper. 

\vspace{3mm}
\noindent
{\it Remark 4.3} In the $\uqslth$ limit, the type I vertex operator 
$\Philz{l}$ \refeq{typeivo} coincides with 
the result in \cite{Matsuo}. 
On the other hand, in the CFT limit, the same
vertex operator coincides with the primary field 
$\Psi_{l;l+1,1}(z)$ in \refeq{primaryf}, whereas
the type II vertex operator $\Psilsz{l}$ coincides with $
\phi_{l,-l}(z):\exp\sqrt{2}\alpha_{1,l+1}\phi_0(z):$. 
Hence one can regard $\Philz{l}$ and $\Psilsz{l}$ as 
the $q-$deformation of the primary fields in the coset theory.

\vspace{3mm}
\noindent
{\it Remark 4.4} At level one, the vertex operator 
$\tilde\Phi^{(1)}_1(z)$ coincides with $\Psi^+(z)$
obtained by Lukyanov and Pugai in the $q-$Virasoro algebra\cite{LukPug}. 
On the other hand, the degeneration limit \refeq{dtypeIt} and 
\refeq{dtypeIIt}, at level one, coincide with those found in 
the massless XXZ model\cite{JKM}. 
In this way, in all the known cases, 
the vertex operators  relevant for the 
physical applications obey the twisted intertwining relations.  
 
The vertex operators \refeq{typeivo}--\refeq{ttypeiivo} satisfy  
interesting commutation relations with the 
currents which allows  us to expect the existence of 
the $(q,p)-$analogue of the ( twisted ) 
intertwining relations \refeq{inthi}-\refeq{tinteii} 
and the existence of the Hopf algebra structure in 
$\uxpslt$ (see also Sec.6).
We here list them only for the vertices $\tilde\Philz{l}$ 
and $\tilde\Psilsz{l}$.

\begin{prop}
The vertex operators $\tilde\Philz{l}$ and $\tilde\Psilsz{l}$
satisfy the following relations.
\bea
&&H^{(\pm)}(w)\tilde\Philz{l}=q^l\Bigl(\frac{z}{w}q^{\pm k/2}\Bigr)^{-l/r}
\frac{\Theta_p(q^{-l\pm k/2}z/w)}{\Theta_p(q^{l\pm k/2}z/w)}\tilde\Philz{l}
H^{(\pm)}(w)
\\
&&E(w)\tilde\Philz{l}+\tilde\Philz{l}E(w)=0,\\
&&F(w)\tilde\Philz{l}=-q^l\Bigl(\frac{z}{w}\Bigr)^{-l/r}
\frac{\Theta_p(q^{-l}z/w)}{\Theta_p(q^{l}z/w)}\tilde\Philz{l}F(w),\\
&&H^{(\pm)}(w)\tilde\Psilsz{l}=q^{-l}\Bigl(\frac{z}{w}q^{\mp k/2}
\Bigr)^{l/(r-k)}
\frac{\Theta_{p^*}(q^{l\mp k/2}z/w)}
{\Theta_{p^*}(q^{-l\mp k/2}z/w)}\tilde\Psilsz{l}
H^{(\pm)}(w),\\
&&F(w)\tilde\Psilsz{l}+\tilde\Psilsz{l}F(w)=0,\\
&&E(w)\tilde\Psilsz{l}=-q^{-l}\Bigl(\frac{z}{w}\Bigr)^{l/(r-k)}
\frac{\Theta_{p^*}(q^{l}z/w)}{\Theta_{p^*}(q^{-l}z/w)}\tilde\Psilsz{l}E(w).
\ena
\end{prop}

Finally, we present the commutation relations among the  type I and type II
vertex operators. For application to physics, we are interested in those
among $\tilde\Phi^{(k)}_k(z)$ and $\tilde\Psi^{(1)}_1(z)$. 
\begin{prop}
\bea
&&\tilde\Phi^{(k)}_{k}(z)\tilde\Phi^{(k)}_k(w)=r_{k}(w/z)
\tilde\Phi^{(k)}_k(w)\tilde\Phi^{(k)}_{k}(z),\\
&&\tilde\Psi^{(1)*}_{1}(z)\tilde\Psi^{(1)*}_1(w)=s_1(z/w)
\tilde\Psi^{(1)*}_1(w)\tilde\Psi^{(1)*}_{1}(z),\\
&&\tilde\Phi^{(k)}_{k}(z)\tilde\Psi^{(1)*}_1(w)=\chi(w/z)
\tilde\Psi^{(1)*}_1(w)\tilde\Phi^{(k)}_{k}(z),
\ena
where
\bea
&&r_k(z)=z^{\frac{k(r-k)}{2r}}
\frac{(q^{2r-2k+2}/z;q^4,p)_\infty\infq4p{q^{2k+2}/z}
(q^{2r+2}z;q^4,p)_\infty\infq4p{q^2z}}
{(q^{2r+2}/z;q^4,p)_\infty\infq4p{q^{2}/z}
(q^{2r-2k+2}z;q^4,p)_\infty\infq4p{q^{2k+2}z}},\nonumber\\
\\
&&s_1(z)=z^{\frac{r}{2k(r-k)}-\frac{(k-1)^2}{k(k+2)}}\frac{\sigma(w/z)}
{\sigma(z/w)}
\ena
with 
\bea
&&\sigma(z)=
\frac{
(q^{2(k+1)}z;q^{2k},q^{2(k+2)})^2_\infty
(q^{2(k+2)}z;q^4,q^{2k})_\infty(q^{2k}z;q^4,q^{2k})_\infty
(q^{4}z;q^4,p^*)_\infty(z;q^4,p^*)_\infty
}
{
(q^{4k}z;q^{2k},q^{2(k+2)})^2_\infty
(q^{4}z;q^{2k},q^{2(k+2)})^2_\infty
(q^{2(k+1)}z;q^4,q^{2k})^2_\infty(q^{2}z;q^4,p^*)^2_\infty
}\nonumber\\
\ena
and
\bea
&&\chi(z)=z^{1/2}\frac{\Theta_{q^4}(q/z)}{\Theta_{q^4}(qz)}.
\ena
\end{prop}
The function $r_k(z)\equiv 1/\kappa(u)$ with $z=q^{2u}$
 satisfies the following inversion relations.
\bea
&&\kappa(u)\kappa(-u)=1,\\
&&\kappa(u)\kappa(-2-u)=\frac{\vartheta_1(\frac{k+1+u}{r})
\vartheta_1(\frac{k-1-u}{r})} {\vartheta_1(\frac{1+u}{r})
\vartheta_1(\frac{-1-u}{r})}.
\ena
Therefore, according to Appendix D in the third reference  
in \cite{DJKMO}, one can identify $\kappa(u)$ with the partition function per 
site for the $k-$fusion RSOS model in the regime III. On the other hand, 
the logarithmic derivative of the function $\chi(z)$ gives the excitation 
energy of the kink\cite{MW}.

\setcounter{equation}{0}
\section{$q-$Deformation of the Coset Theory}
In this section, we discuss a $q-$deformation of the coset conformal field 
theory  based on the algebra $\uxpslt$ and make an identification of it 
with the $k-$fusion RSOS model $k\in\Z_{>0}$.

\subsection{Definition}
The Virasoro central charge of the $\Z_k$ parafermion theory is 
$c_{PF}=\frac{2(k-1)}{k+2}$.
Hence the level one $(k=1)$  parafermion theory has  zero 
central charge and gives a trivial contribution to the coset theory.  
This should be true in the $q-$deformed theory, too. 
Noting this, one can see that at level one, the currents $E(z)$ and $F(z)$ 
in \refeq{curre} and \refeq{currf} 
coincides with the screening currents $\bar{V}(z)$
and $\bar{U}(z)$ in the $q-$Virasoro algebra\cite{LukPug}.
In addition, the type I vertex $\tilde\Philz{1}$ 
coincides with the vertex  $\Psi^+(z)$ in \cite{LukPug,AKOS}.  
In this sense, the algebra $\uxpslt$ at level one governs 
the structure of the $q-$Virasoro algebra. 
In fact, at level one, the $q-$Virasoro algebra generator $T(z)$
is obtained by\cite{AKOS}
\bea
&&T(z)=\Lambda(zq)+\Lambda(zq^{-1}),\\
&&\Lambda(z)=:\tilde\Phi^{(1)}_1(zq)\tilde\Phi^{(1)}_1(zq^{-1}):.
\ena

For generic level $k$, as mentioned in Remarks 4.1 and 4.3, one can 
regard the currents $E(z)$, $F(z)$ and the type I vertex $\Philz{l} 
(l=0,1,..,k)$
as the $q-$deformation of the screening currents $S_+(z)$, 
$S_-(z)$ and the primary field $\Psi_{l;l+1,1}$ in the coset CFT. 

These observation lead us to the 
following characterization of a  $q-$deformation of 
the coset theory $SU(2)_k\times SU(2)_{r-k-2}/SU(2)_{r-2}$,
or equivalently a $q-$deformation of the extended Virasoro algebra
in the free boson realization. 
\begin{enumerate}
\item
The theory is obtained as a composition of the $q-$deformed $\Z_k$
parafermion theory and the $\phi_0$ boson theory.
\item The screening currents  satisfy the algebra 
$\uxpslt$.
\item 
The $q-$deformation of the primary fields 
are the intertwiners between the infinite dimensional 
representations of $\uxpslt$. 
\end{enumerate}

We have not yet succeeded to obtain a $q-$Virasoro
generators with central charge $c$ \refeq{centralc} and any extra generators. However, the free boson realization of 
the screening operators and 
the type I,  type II vertex operators 
enables us to  analyze the 
structure of the highest weight representation of 
the  $q-$Virasoro algebras\cite{LukPug}.
In the following section, 
we carry out  such analysis for the case $k>1$. The resultant irreducible representation 
turns out to be identified with the space of states 
of the $k-$fusion RSOS model. 

\subsection{Fock modules }

Let $J=|n'-n\ ({\rm mod} 2k)|$ and $M=n'-n$ (mod $2k$). 
Define the highest weight state $\ket{J,M;n',n}$ by
\bea
&&\ket{J,M;n',n}=\ket{J,M}_{PF}\otimes \ket{n',n}_0,\\
&&\ket{J,M}_{PF}=e^{\frac{J}{2(k+2)}Q_1+\frac{M}{2k}Q_2}\ket{0}_{PF},\\
&&\ket{n',n}=e^{-i\sqrt{2}\alpha_{n',n}Q_0}\ket{0}_0. 
\ena
Here $\ket{0}_{PF}$ and $\ket{ 0}$ denote the $SL(2,\C)$ invariant 
vacuum states defined by
\bea
&& a_{j,m}\ket{0}_{PF}=0=P_j\ket{0}_{PF} \qquad (j=1,2\ m\in\Z_{>0}),\\
&& a_{0,m}\ket{0}_{0}=0=P_0\ket{0}_{0}\qquad ( m\in\Z_{>0}) .
\ena
Note that the highest weight states can be obtained 
by making the vertex operator act on the $SL(2,\C)$ invariant vacuum state
in the same way as \refeq{hwstate} in CFT.

The  Fock modules $\F_{J,M;n',n}=\F^{PF}_{J,M}\otimes 
\F^{\phi_0}_{n',n}$ are defined by
\bea
&&\F^{PF}_{J,M}=\C[a_{1,-m_1},\ a_{2,-m_2}\ (m_1,m_2\in\Z_{>0})]\ket{J,M}_{PF}
,\\
&&\F^{\phi_0}_{n',n }=\C[a_{0,-m}\ (m \in\Z_{>0})]\ket{n',n}_{0}.
\ena

We also define the degree of the vector in the Fock modules 
as an eigenvalue of the  operator $L_0$ given by
\bea
&&L_0=L_0^{PF}+L_0^{\phi_0},\label{lzero}\\
&&L^{PF}_0=\sum_{m>0}\frac{m^2}{[2m][(k+2)m]}a_{1,-m}a_{1,m}+
\frac{P_1(P_1+2)}{4(k+2)}\nonumber\\
&&\qquad\qquad-\sum_{m>0}\frac{m^2}{[2m][km]}a_{2,-m}a_{2,m}-
\frac{P_2^2}{4k},\\
&&L_0^{\phi_0}=\sum_{m>0}\frac{m^2[(r-k)m]}{[2m][km][rm]}a_{0,-m}a_{0,m}+
\frac{1}{2}{P_0\Bigl(P_0-\sqrt{\frac{2k}{r(r-k)}}\Bigr)}.
\ena
For a vector $u\in \F_{J,M;n',n} $
\bea
&& u=\biggl(
\prod_{i=0,1,2}a_{i,-m_{i,1}}a_{i,-m_{i,2}}\cdots a_{i,-m_{i,N_i}}
\biggr)\ket{J,M;n',n}, 
\ena
its degree $N$ is given by
\bea  
&&L_0 u=(\ h_{J,M;n',n}+N )u, \qquad 
N=\sum_{i=0,1,2}\sum_{j=1}^{N_i} m_{i,j}
\label{degreecount}\ena
where
\bea
&&h_{J,M;n',n}=h_{J,M}+
{(nr-n'(r-k))^2-k^2\over 4kr(r-k)},\label{hcosetjm}\\
&&h_{J,M}={J(J+2)\over 4(k+2)}-\frac{M^2}{4k}\label{hjm}.
\ena

\subsection{Screening currents}
In section 5.1, we identified the currents  $E(z)$ and $F(z)$ 
with the $q-$deformation of the screening currents $S_+(z)$
and $S_-(z)$ in the coset CFT. As mentioned in section 2, we need two more 
screening currents, which governs the structure of the Fock 
representation space of the $q-$deformed $\Z_k$ parafermion
theory. Such  screening currents $S(z):\F_{J,M;n,n'}\to\F_{J-2,M;n,n'}$ and 
$\eta(z):\F_{J,M;n,n'}\to\F_{J+k+2,M+k;n,n'}$ 
were already realized by Matsuo\cite{Matsuo}.
They are given by the following formulae.
\bea
&&S(z)=\frac{-1}{(q-q^{-1})}
:\EXP{\phi_1\Bigl(k+2\Bigl|z;-\frac{k+2}{2}\Bigr)}\nonumber\\
&&\qquad\qquad\times\Bigl(
\EXP{\phi_2^{(+)}\Bigl(1;2\Bigl|z;\frac{k+2}{2}\Bigr)
+\phi_1^{(+)}\Bigl(1;2\Bigl|z;\frac{k}{2}\Bigr)}\nonumber\\
&&\qquad\qquad\qquad\qquad-\EXP{-\phi_2^{(-)}\Bigl(1;2\Bigl|z;\frac{k+2}{2}\Bigr)
-\phi_1^{(-)}\Bigl(1;2\Bigl|z;\frac{k}{2}\Bigr)}\Bigr):,
\\
&&\eta(z)=:\EXP{-\phi_1\Bigl(2\Bigr|z;\frac{k}{2}\Bigr)
-\phi_2\Bigl(2\Bigr|z;\frac{k+2}{2}\Bigr)}
:.
\ena

Noting that $S(z)$ and $\eta(z)$ depend only on the boson fields 
$\phi_1$ and $\phi_2$, the following relations are direct consequences  of the
Lemma 4.1 and 4.5 in \cite{Matsuo}.
\begin{prop}
\bea
&&E(z)S(w)=S(w)E(z)=O(1),\\
&&F(z)S(w)=S(w)F(z)=[k+2]{}_{k+2}\partial_{w}\Bigl[
\frac{1}{z-w}:e^{\Phi_{FS}(z)}:\Bigr]+O(1)\\
&&E(z)\eta(w)=-\eta(w)E(z)=
{}_{1}\partial_{w}\Bigl[\frac{1}{z-w}:e^{\Phi_{E\eta}(z)}:\Bigr]+O(1),\\
&&F(z)\eta(w)=-\eta(w)F(z)=O(1),\\
&&S(z)\eta(w)=\eta(w)S(z)=
{}_{1}\partial_{w}\Bigl[
\frac{1}{z-w}:e^{\Phi_{S\eta}(z)}:\Bigr]+O(1)
\ena
with
\begin{eqnarray*}
&&\Phi_{FS}(z)=\boszerop{k}{z}+
\bosone{k+2}{z}{\frac{k+2}{2}}+\bostwo{k}{z}{\frac{k}{2}},\\
&&\Phi_{E\eta}(z)=-\boszero{k}{z}-
\bosone{2}{z}{\frac{k+2}{2}}-\bostwo{k}{z}{-\frac{k}{2}}
-\bostwo{2}{z}{\frac{k+4}{2}},\\
&&\Phi_{S\eta}(z)=
\bosone{k+2}{z}{-\frac{k+2}{2}}-\bosone{2}{z}{\frac{k-2}{2}}
-\bostwo{2}{z}{\frac{k}{2}}.
\end{eqnarray*}
Here the difference of the function $f(z)$ is defined by
$$
{}_a\partial_z \ f(z)=\frac{f(q^a z)-f(q^{-a} z)}{q-q^{-1}}.
$$
\end{prop}

In addition, from Theorem 4.1, the following commutation relations hold. 
\begin{prop}
\bea
&&E(z)E(w)=\zfactor{u-v+1}{u-v-1}E(w)E(z),\\
&&F(z)F(w)=\zfactorp{u-v-1}{u-v+1}F(w)F(z),\\
&&S(z)S(w)=\zfactorpf{u-v+1}{u-v-1}S(w)S(z).
\ena
\end{prop}
Here the notation $[u]_x \ (x=r, r-k, k+2)$ stands for the theta function
\bea
&&[u]_x=\vartheta_1\Bigl(\ \frac{u}{x}\ \Bigl|\ \tau_x\ \Bigr),
\ena
where we set $q^{2x}=e^{-2\pi i/\tau_x}$.
Hence $\tau_r=\tau$ and $ \tau_{r-k}=\tau^*$.
One should note the following quasi-periodicities.
\bea
&&[u+x]_x=-[u]_x,\qquad [u+x\tau_{x}]_x=-e^{-\pi i(2u+\tau_x)}[u]_x.
\label{quasiperiod}
\ena

Now let us define a set of screening operators. 
\begin{dfn}
\bea
&&Q_+=\Ointz{z}E(z)\zfactor{u-\frac{1}{2}+\hat\Pi}{u-\frac{1}{2}},\\
&&Q_-=\Ointz{z}F(z)\zfactorp{u+\frac{1}{2}-\hat\Pi'}{u+\frac{1}{2}},\\
&&Q=\Ointz{z}S(z)\zfactorpf{u-\frac{1}{2}+P_1}{u-\frac{1}{2}},
\label{screenqs}\\
&&\eta_0=\Oint{z} \eta(z),
\ena
where 
\bea
&&\hat\Pi=\sqrt{\frac{2r(r-k)}{k}}P_0-\frac{r-k}{k}P_2,\\
&&\hat\Pi'=\sqrt{\frac{2r(r-k)}{k}}P_0-\frac{r}{k}P_2.
\ena
\end{dfn}

Due to the quasi-periodicity \refeq{quasiperiod}, the integrands 
in $Q_+,\ Q_-$ and $ Q$ are single valued in $z$. In addition
the integrand in $\eta_0$ is single valued on $\F_{J,M;n',n}$. 
Therefore all the integrations in Definition 5.3 
can be  taken over a closed contour on $\F_{J,M;n',n}$.

The following commutation relations hold.
\begin{lem}
\bea
&& [Q_\pm,Q]=0,\label{qpmq}\\
&& [Q_\pm,\eta_0]=0,\\
&& \{Q,\eta_0\}=0,
\ena
where $\{A,B\}=AB+BA$.
\end{lem}

\noindent
{\it Proof of \refeq{qpmq}}
From Proposition 5.1 and $[\hat\Pi,S(w)]=0,\ [P_1,E(z)]=0$, 
the relation $[Q_+,Q]=0$ is trivial.
For the commutation $[Q_-,Q]$, we have from Proposition 5.1 and 
 $[\hat\Pi',S(w)]=0,\ [P_1,F(z)]=0$
\bea
&&[Q_-,Q]=-[k+2]\OintC{|z|=1}{z}\OintC{C_z}{w} 
{}_{k+2}\partial{w}\Bigl[
\frac{1}{z-w}:e^{\Phi_{FS}(z)}:\Bigr]
\nonumber\\
&&\qquad\qquad\times\zfactorp{u+\frac{1}{2}-\hat\Pi'}{u+\frac{1}{2}}
\zfactorpf{v-\frac{1}{2}+\Pi_1}{v-\frac{1}{2}}.
\label{proofqpmq}
\ena
Here $C_z$ denotes a closed contour enclosing the points $q^{\pm (k+2)}z$. 
After taking the integral 
in $w$, the quasi-periodicity of the theta function $[u]_{k+2}$
makes the right hand side of \refeq{proofqpmq} vanish.  

The other statements can be proved in the similar way.
\qed

\vspace{3mm}
\noindent
{\it Remark 5.1} Due to the theta function factor in the Definition 5.3,
the screening operators $Q_+$ and $Q_-$ does not commute each other. 
However, the relation $[Q_+^n,Q_-^{n'}]=0$ holds on $\F_{J,M;n',n}$. 

Let us set $Q^+_n=Q^n_+,\ Q^-_n=Q^n_-,$ and $Q_n=Q^n$. Then we claim
\begin{thm}
The screening operators $Q_+,\ Q_-,\ Q$ and $\eta_0$ are nilpotent:
\bea
&&Q^+_nQ^+_{r-k-n}=Q^+_{r-k-n}Q^+_n=Q^+_{r-k}=0,\\
&&Q^-_nQ^-_{r-n}=Q^-_{r-n}Q^-_n=Q^-_{r}=0,\\
&&Q_nQ_{k+2-n}=Q_{k+2-n}Q_n=Q_{k+2}=0,\\
&&\eta_0^2=0.
\ena
\end{thm}
The proof is due to the following lemma\footnote{
The screening operator $S(z)$ is equivalent to the  one 
in $\uqslth$ discussed in \cite{Konno94}.  We hence have proved the nilpotency of the BRST charge in $\uqslth$ in the improved  form as \refeq{screenqs}. }.
\begin{lem}
\bea
&&Q_n^{+}=\frac{[n]_{r-k}!}{n![1]_{r-k}^n}\prod_{j=1}^n\biggl(
\Ointz{z_j}
E(z_j)\biggr)\prod_{i<j}\frac{[u_i-u_j]_{r-k}}{[u_i-u_j-1]_{r-k}}
\prod_{i=1}^n
\frac{[u_i+\frac{1}{2}+\hat\Pi-n]_{r-k}}{[u_i-\frac{1}{2}]_{r-k}},\nonumber\\
\\
&&Q_n^{-}=\frac{[n]_r!}{n![1]_r^n}
\prod_{j=1}^n\biggl(
\Ointz{z_j}
F(z_j)\biggr)
\prod_{i<j}\frac{[u_i-u_j]_r}{[u_i-u_j-1]_r}\prod_{i=1}^n
\frac{[u_i-\frac{1}{2}-\hat\Pi'+n]_r}{[u_i+\frac{1}{2}]_r}, \\
&&Q_n^{}=\frac{[n]_{k+2}!}{n![1]_{k+2}^n}
\prod_{j=1}^n\biggl(
\Ointz{z_j}
S(z_j)\biggr)
\prod_{i<j}\frac{[u_i-u_j]_{k+2}}{[u_i-u_j-1]_{k+2}}
\prod_{i=1}^n
\frac{[u_i+\frac{1}{2}+P_1-n]_{k+2}}{[u_i-\frac{1}{2}]_{k+2}}.\nonumber\\
\ena
\end{lem}
This lemma is proved by using the commutation relations in Proposition 5.2 and 
the following theta function identity\cite{JLMP}.
\begin{lem}
\bea
&&\frac{1}{n!}\sum_{\sigma\in S_n}\prod_{i=1}^n[u_{\sigma(i)}-2i+2]_x
\prod_{i<j\atop \sigma(i)>\sigma(j)}
\zfactora{u_{\sigma(i)}-u_{\sigma(j)}-1}
{u_{\sigma(i)}-u_{\sigma(j)}+1}\nonumber\\
&&=\frac{[n]_x!}{n![1]^n_x}\prod_{i<j}
\zfactora{u_{i}-u_{j}}
{u_{i}-u_{j}-1}\prod_{i=1}^n[u_i-n+1]_x \qquad (x=r,\ r-k,\ k+2).
\ena
\end{lem}

\subsection{Resolution of the Fock modules $\F_{J,M;n',n}$}
The Fock modules $\F_{J,M;n',n}=\F^{PF}_{J,M}\otimes \F^{\phi_0}_{n',n}$ 
are reducible due to the existence of the 
singular vectors, which are constructed by the screening operators on 
some highest weight states\cite{Konno94,AKOS}. 
In order to obtain irreducible 
spaces, we consider a resolution of the modules $\F_{J,M;n',n}$
following the method by Felder\cite{Felder}.

As mentioned in Sec.2, Lemma 5.4 allows us to divide the consideration
into the following two steps. First consider the resolution
of the $\Z_k$ parafermion Fock modules $\F^{PF}_{J,M}$  and 
get spaces $\H^{PF}_{J,M}$ as  irreducible representations.
Then consider the resolution of the modules $\F^{PF}_{J,M}\otimes
\F^{\phi_0}_{n',n}$.

\vspace{3mm}
\noindent
{\bf 1) Resolution of the $\Z_k$ $q-$parafermion Fock modules $\F^{PF}_{J,M}$}
\vspace{3mm}

Let us first remind the reader the following two facts. 
\begin{enumerate}
\item The $\Z_k$ parafermion theory is obtained as the coset $SU(2)_k/U(1)$.
In the $q-$deformed case, especially in our case, 
this means that the $q-$deformed $\Z_k$ parafermion theory 
is obtained from Matsuo's bosonization of $\uqslth$
by dropping $\{\ \a_m \}$ boson. 
\item The screening currents $S(z)$ and $\eta(z)$
are independent of the $\{ \a_m\}$ bosons. 
Therefore these screening currents 
are common in the $\uqslth$ theory and the $q-$parafermion theory 
so that  
the structure of the singular vectors in 
the $q-$parafermion Fock modules  are the same as corresponding Fock modules 
for $\uqslth$. 
\end{enumerate}

The structure of the Fock modules for 
$\uqslth$ was investigated in \cite{Matsuo,Konno94}. 
The modules are given by $\F_J=\bigoplus_{M\in\Z}
\F^{PF}_{J,M}\otimes 
\F^{\a}_M$. Here $\F^{\a}_M$ denotes the Fock module 
of the $\{\ \a_m\}$ boson. 
In \cite{Konno94}, we showed, using the 
equivalent boson realization, that the Fock modules $\F_J$ 
are reducible for the highest weight $\lambda_{a,a'}=
J_{a,a'}\Lambda_1+(k-J_{a,a'})\Lambda_0$ due to 
the existence of the singular vectors. 
Here $J_{a,a'}=a-(k+2)a'-1$, 
the level $k$ being , in general, a rational number
$k+2=P/P',\ P,P'\geq 1, \ {\rm GCD}(P,P')=1 $ and
$ 1\leq a\leq P-1, \ 0 \leq a'\leq P'-1$.
The case $P'=1$ is relevant for the parafermion theory.
It was then observed that the Fock module structure, i.e. the 
degree of the singular vectors and 
the cosingular vectors as well as 
the multiplicities in each degree seems to be 
the same as the CFT case\cite{BerFel}. This was checked by calculating the 
characters of the irreducible 
representation spaces obtained 
by a resolution of the Fock modules. 

Following the procedure given for  $\uqslth$\cite{Matsuo,Konno94} 
and the analysis in 
CFT\cite{DistQiu,FLMSS}, 
we make the resolution of the $q-$parafermion Fock modules
$\F^{PF}_{J,M}$ as follows.

We first consider the following restriction.
\bea
&&\tilde\F^{PF}_{J,M}=\Ker(\eta :\ 
\F^{PF}_{J,M}{\to} \F^{PF}_{J+k+2,M+k}).
\ena
Since $\eta^2_0=0$, the complex of the 
Fock modules associated with the map $\eta_0 :
\ \F^{PF}_{J,M}{\to} \F^{PF}_{J+k+2,M+k}$ has trivial cohomology. 
We hence have
\bea 
&&\tr_{\tilde \F^{PF}_{J,M}}{\cal O}=\sum_{u\geq 0}(-)^u
\tr_{\F^{PF[u]}_{J,M}}{\cal O}
\label{vtrace}
\ena
and 
\bea
&&0=\sum_{u\in \bf Z}(-)^u
\tr_{\F^{PF[u]}_{J,M}}\cal O,
\label{vtraceb}
\ena
where
$\F^{PF[u ]}_{J,M}\equiv
\F^{PF}_{J+(k+2)u,M+ku}$.

On the other hand, the 
operator $Q$ generates the following  complex of the 
restricted  Fock
modules $\tilde\F^{PF}_{J,M}$.
\bea
&&\cdots\stackrel{Q^{[-2]}_{J+1}}{\longrightarrow}
\tilde\F^{PF[-1]}_{J,M}\stackrel{Q^{[-1]}_{J+1}}{\longrightarrow}
\tilde\F^{PF[0]}_{J,M}
\stackrel{Q^{[0]}_{J+1}}{\longrightarrow}
\tilde\F^{PF[1]}_{J,M}\stackrel{Q^{[1]}_{J+1}}{\longrightarrow}
\tilde\F^{PF[2]}_{J,M}\stackrel{Q^{[2]}_{J+1}}{\longrightarrow}\cdots.
\label{complexq}\ena
Here we introduced the notations $Q_{J+1}^{[2s]}=Q_{J+1}$ and 
$Q_{J+1}^{[2s+1]}=Q_{(k+2)-J-1}$. These are operators 
acting  on the modules
$\F^{PF[2s]}_{J,M}=\F^{PF}_{J-2s(k+2),M}$ and 
$\F^{PF[2s+1]}_{J,M}=\F^{PF}_{-J-2-2s(k+2),M}$,
respectively.

The  results in CFT\cite{BerFel,DistQiu,FLMSS,Konno92} and 
the investigation in the quantum affine algebra 
$\uqslth$\cite{Konno94} lead us to the following conjecture 
\begin{conj}
The cohomology of the complex \refeq{complexq} is given by
\bea
&&\Ker Q^{[s]}_{J+1}/\Image Q^{[s-1]}_{J+1}=
\Biggl\lbrace\matrix{0&&for\quad s\not=0\cr
                      {\cal H}^{PF}_{J,M}&&for \quad s=0 \cr },
\label{paracohom}
\ena
where ${\cal H}^{PF}_{J,M}$ is the irreducible highest weight 
module of the $q-$deformed $\Z_k$ parafermion theory  with the highest weight 
$h_{J,M}$ \refeq{hjm}.

\end{conj}

As a consequence, we obtain the following trace formula
\begin{cor}
\bea
&&\tr_{{\cal H}^{PF}_{J,M}}
{\cal O}=\sum_{s\in \Z}(-)^s
\tr_{\tilde \F^{PF[s]}_{J,M}}
 {\cal O}^{\lbrack s\rbrack}, \label{paratrace}
\ena
where ${\cal O}^{\lbrack s\rbrack}$ 
is an operator on $\F^{PF[s]}_{J,M}$ obtained
by the recursion formula
\bea
&&Q^{\lbrack s\rbrack}_{a} {\cal O}^{\lbrack s\rbrack}=
{\cal O}^{\lbrack s+1\rbrack}
Q^{\lbrack s\rbrack}_{a} 
\ena
with ${\cal O}^{\lbrack 0\rbrack}={\cal O}$.
\end{cor}

Combining \refeq{paratrace} and \refeq{vtrace}, we get 
\bea
&&
\tr_{{\cal H}^{PF}_{J,M}}{\cal O}
=\sum_{s\in \bf Z}\sum_{u\geq 0}(-)^{s+u}
\tr_{ \F^{PF[s,u]}_{J,M}} 
{\cal O}^{\lbrack s\rbrack},
\label{cparatrace}
\ena
where 
\bean
&&\F^{PF[2s,u]}_{J,M}=\F^{PF}_{J+(k+2)(u-2s),M+ku},\\ 
&&\F^{PF[2s+1,u]}_{J,M}=\F^{PF}_{-J-2+(k+2)(u-2s),M+ku}.
\enan 

To support the conjecture, let us apply this formula to the calculation of  the character of the space ${\cal H}^{PF}_{J,M}$.
We obtain
\bea
&&\chi_{J,M}(\omega)=\tr_{{\cal H}^{PF}_{J,M}}\ y^{L^{PF}_0-c_{PF}/24}
=\eta(\omega)c^{J}_{M}(\omega) \label{pfchar}
\ena
with $y=e^{2\pi i\omega}$.
Here the function $c^J_{M}(\omega)$ denotes  the string function 
\bea
&&c^J_{M}(\omega)=\eta(\omega)^{-3}\sum_{s\in\Z}\sum_{u\geq 0}(-1)^u
\Bigl(x^{B^{[s,u]}_{J,M}}-x^{B^{[s,u]}_{-J-2,M}}\Bigr),\\
&&B^{[s,u]}_{J,M}=\frac{(J+1-2(k+2)(s-u/2))^2}{4(k+2)}-\frac{(M+ku)^2}{4k}
\label{stringf}
\ena
with  $\eta(\omega)$ being
the Dedekind eta function
\bea
&&\eta(\omega)=y^{1/24}\prod_{n=1}^\infty(1-y^n).
\ena
The result \refeq{pfchar} is precisely the irreducible character of 
the $\Z_k$ parafermion theory.

\vspace{3mm}
\noindent
{\bf 2) Resolution of the Fock modules $\F_{J,M;n',n}$}
\vspace{3mm}

Replacing the space $\F^{PF}_{J,M}$ with ${\cal H}^{PF}_{J,M}$ in
$\F_{J,M;n',n}$,
now we consider the Fock modules 
$\tilde\F_{J,M;n',n}={\cal H}^{PF}_{J,M}\otimes \F^{\phi_0}_{n',n}$
with $J=|n'-n\ ({\rm mod}2k)|$, $M=n'-n$ ( mod $2k $),
 $1\leq n\leq r-k-1,\ 1\leq n'\leq r-1, 
0\leq J\leq k$.

Let us introduce the notations
 $Q^{+[2t]}_n=Q^+_n$, $Q^{+[2t+1]}_n=Q^+_{r-k-n}$, 
 $Q^{-[2t]}_{n'}=Q^-_{n'}$, $Q^{-[2t+1]}_{n'}=Q^-_{r-n'} \ (t\in\Z)$ and 
\bea
&&\tilde\F^{[2t]}_{J,M;n',n}=\tilde\F^{}_{J,M+2(r-k)t\ ;n',n-2(r-k)t},\\
&&\tilde\F^{[2t+1]}_{J,M;n',n}=
\tilde\F^{}_{J,M+2n+2(r-k)t\ ;n',-n-2(r-k)t}\\
&&\tilde\F^{[2t]'}_{J,M;n',n}=\tilde\F^{}_{J,M+2 rt\ ;n'-2rt,n},\\
&&\tilde\F^{[2t+1]'}_{J,M;n',n}=
\tilde\F^{}_{J,M-2n'+2rt\ ;-n'-2rt, n}.
\ena
Then,  
the screening operators $Q_+, \ Q_-$ generate the following 
infinite sequences of the Fock
modules $\tilde\F_{J,M;n',n}$.
\bea
&&\cdots\stackrel{Q^{+[-2]}_{n}}{\longrightarrow}
\tilde\F^{[-1]}_{J,M;n',n}\stackrel{Q^{+[-1]}_{n}}{\longrightarrow}
\tilde\F^{[0]}_{J,M;n',n}
\stackrel{Q^{+[0]}_{n}}{\longrightarrow}
\tilde\F^{[1]}_{J,M;n',n}\stackrel{Q^{+[1]}_{n}}{\longrightarrow}\cdots,
\label{complexp}\\
&&\cdots\stackrel{Q^{-[-2]}_{n'}}{\longrightarrow}
\tilde\F^{[-1]'}_{J,M;n',n}\stackrel{Q^{-[-1]}_{n'}}{\longrightarrow}
\tilde\F^{[0]'}_{J,M;n',n}
\stackrel{Q^{-[0]}_{n'}}{\longrightarrow}
\tilde\F^{[1]'}_{J,M;n',n}\stackrel{Q^{-[1]}_{n'}}{\longrightarrow}\cdots.
\label{complexm}
\ena
Due to Theorem 5.5, 
these  are complexes.
As in the level one case\cite{LukPug},
it is enough to consider one of them. Let us consider the complex
\refeq{complexp}.
The following observation suggests  the existence of the singular vectors 
in $\tilde\F_{J,M;n',n}$
similar to the CFT case\cite{Konno92}.

Let  $M_{n',n}=n'-n$ (mod $2k$) and consider the vector
\bea
&&\ket{\chi_{-n',n}}=Q_n^+\ket{J,M_{-n',n};-n',n}\in 
\F_{J,M_{-n',-n};-n',-n}.
\ena
Using Lemma 5.6 and the operator product
\bea
&&\Psi(z_1)\Psi(z_2)
=\Bigl(\frac{1}{q-q^{-1}}\Bigr)^2 z_1^{-\frac{2}{k}}q
\frac{(q^{2(k+1)}z_2/z_1;q^{2k})}{(q^{-2}z_2/z_1;q^{2k})}\nonumber\\
&&\qquad\times\Bigl(
(1-z_2/z_1):\Psi_I(z_1)\Psi_{I}(z_2):-(1-q^{-2}z_2/z_1)
:\Psi_I(z_1)\Psi_{II}(z_2):
\nonumber\\
&&\qquad-q^{-2}(1-q^{2}z_2/z_1):\Psi_{II}(z_1)\Psi_{I}(z_2):
+(1-z_2/z_1):\Psi_{II}(z_1)\Psi_{II}(z_2):\Bigr),
\ena
where we set the RHS of \refeq{pfcurrents} as 
$\frac{1}{q-q^{-1}}(\Psi_I(z)-\Psi_{II}(z))$,
we have 
\bea
&&Q^+_n\ket{J,M_{-n',n};-n',n}\nonumber\\
&&=\frac{[n]_{r-k}!}{n![1]_{r-k}^n}\Bigl(\frac{-1}{q-q^{-1}}\Bigr)^n
\Ointz{z_1}\cdots\Ointz{z_n}
\prod_{i=1}^n z_i^{2\sqrt{\frac{r}{k(r-k)}}\a_{n',n}-\frac{M}{k}}
\prod_{i<j}z_i^{\frac{2}{r-k}}q
\frac{(q^{2(k+2)} z_j/z_i;q^{2k})}{(q^{-2}z_j/z_i;q^{2k})}
\nonumber\\
&&\times\ {\rm sum\ of }\biggl[({\rm  polynomials\ of}\ z_j/z_i\ (1\leq i<j\leq n) )
\nonumber\\
&&\times\EXP{{\rm polynomials\ of}\ a_{0,-m}, a_{1,-m'}, a_{2,-m''}, 
z_j\ ( m,m',m''\in\Z>0,\ j=1,..,n )}
\biggr]\nonumber\\
&&\qquad\qquad\qquad \times\ket{J,M_{-n',-n};-n',-n}.\label{qket}  
\ena
Setting $z_1=z,\ z_j=z w_j\ (j=2,..,n)$ and collecting all the $z-$dependence
in the integrand, one can factor the integral   
$\oint dz\ z^{-1-N-\frac{n(n'+n)}{k}-\frac{n}{k}M_{-n',n}}$, where  
$N\in\Z_{\geq 0}$ comes from the exponent in the fourth line of \refeq{qket}. 
Let us evaluate $N$. From \refeq{degreecount}, the non-vanishing term in 
\refeq{qket} has degree $N+h_{J,M_{-n',-n};-n',-n}$. On the other hand,
the BRST charge $Q^+_n$ commutes with $L_0$. Hence the degree of the same term
should be equal to $h_{J,M_{-n',n};-n',n}$. Hence we have
$N=h_{J,M_{-n',n};-n',n}-h_{J,M_{-n',-n};-n',-n}=
(4n'n+M_{-n',-n}^2-M_{-n',n}^2)/4k$. This value $N$ is consistent with
the non-vanishment of the integral $\oint dz$.

We hence  conjecture the following statement\cite{Konno92}. 

\begin{conj}
The complex \refeq{complexp}
 has the non-trivial cohomology only at 
$t=0$, i.e.
\bea
&&\Ker Q_n^{+[t]}/ \Image Q_n^{+[t-1]}
=\Biggl\lbrace\matrix{0&\qquad\qquad for\quad t\not=0\cr
{\cal H}_{J;n',n}&\qquad\qquad for\quad t=0\cr}.
\label{coscohom}
\ena
Here the space ${\cal H}_{J;n',n}$ is the conjectural
irreducible highest weight representation space of the
$q-$Virasoro algebra with central charge $c$ \refeq{centralc} with
the highest weight $h_{J;n',n}=h_{J,J;n',n}$.
\end{conj}
 
As a consequence,  we can derive a trace formula, which 
relates the trace over ${\cal H}_{J;n',n}$ to those over the 
Fock spaces $\tilde\F^{[t]}_{J,M;n',n}\ t\in\Z$. 

\begin{cor}
\bea
&&\tr_{{\cal H}_{J;n,n'}}{\cal O}=
\sum_{t\in \bf Z}(-)^{t}
\tr_{\tilde\F^{[t]}_{J,M;n',n}} 
{\cal O}^{\lbrack t\rbrack}
, \label{costraceb}
\ena
where operator ${\cal O}^{\lbrack t\rbrack}$ on $\tilde \F^{[t]}_{J,M;n',n}$ 
is defined recursively by
\bea
&&Q^{+\lbrack t\rbrack}_n {\cal O}^{\lbrack t\rbrack}=
{\cal O}^{\lbrack t+1\rbrack}
Q^{+\lbrack t\rbrack}_n \label{cosrec}
\ena
with ${\cal O}^{\lbrack 0\rbrack}={\cal O}$.
\par
\end{cor}

      Combining \refeq{cparatrace} and \refeq{costraceb}, we finally get 
the formula
\bea
&&\tr_{{\cal H}_{J;n',n}}{\cal O}=\sum_{s\in \Z}
\sum_{t\in \bf Z}\sum_{u\geq 0}(-)^{s+t+u}
\tr_{\F^{[s,t,u]}_{J,M}\otimes \F^{\phi_0[t]}_{n',n}} 
{\cal O}^{\lbrack s,t\rbrack}
\label{cosettrace}
\ena
where
the  Fock modules $\F^{[s,t,u]}_{J,M;n',n}$ denote the following.
\bea
&&\F^{[2s,2t,u]}_{J,M;n',n}=\F^{PF[u]}_{J-2s(k+2),M+2t(r-k)}\otimes\F^{\phi_0}_{n',n-2(r-k)t},\\
&&\F^{[2s+1,2t,u]}_{J,M;n',n}=
\F^{PF[u]}_{-J-2-2s(k+2),M+2t(r-k)}\otimes\F^{\phi_0}_{n',n-2(r-k)t},
\label{deffock}\\
&&\F^{[2s,2t+1,u]}_{J,M;n',n}=
\F^{PF[u]}_{J-2s(k+2),M+2n+2t(r-k)}\otimes\F^{\phi_0}_{n',-n-2(r-k)t},\\
&&\F^{[2s+1,2t+1,u]}_{J,M;n',n}=
\F^{PF[u]}_{-J-2-2s(k+2),M+2n+2t(r-k)}\otimes\F^{\phi_0}_{n',-n-2(r-k)t}.
\ena
The  operator  
${\cal O}^{[s,t]}$ acting on the space $\F^{[s,t,u]}_{J,M;n',n}$
is defined recursively by
\bea
&&{\cal O}^{[s,0]}={\cal O}^{[s]},\label{obdef}\\
&&Q^{+\lbrack t\rbrack}_n {\cal O}^{\lbrack s,t\rbrack}=
{\cal O}^{\lbrack s,t+1\rbrack}
Q^{+\lbrack t\rbrack}_n. 
\label{obdefb}
\ena


Applying the formula \refeq{cosettrace}, we obtain
 the 
character of the  space ${\cal H}^{}_{J;n',n}$:
\bea
&&\chi_{J;n',n}(\omega)=\tr_{{\cal H}_{J;n,n'}}
y^{L_0-c/24}\nonumber\\
&&\qquad\qquad=\sum_{\bar{M}=-k+1}^k\ c_{J,\bar{M}}
(\omega)\Bigl( \sum_{t\in\bf Z}\delta_{\bar{M},M_{n',n-2(r-k)t}}
y^{B_{\phi_0}^{[2t]}}-\sum_{t\in\bf Z}\delta_{\bar{M},M_{n',-n-2(r-k)t}}
y^{B_{\phi_0}^{[2t+1]}}\Bigr),\nonumber\\
\label{charact}
\ena
where $c$ is given by \refeq{centralc} and 
\bea
&&B_{\phi_0}^{[2t]}=\frac{1}{4kr(r-k)}\Bigl(n'(r-k)-nr+2r(r-k)t\Bigr)^2,\\
&&B_{\phi_0}^{[2t+1]}=\frac{1}{4kr(r-k)}\Bigl(n'(r-k)+nr+2r(r-k)t\Bigr)^2.
\ena
The result \refeq{charact} precisely gives the branching 
coefficient in the formula 
\bea
\chi_{J}^{(k)}(\omega)\chi^{(l)}_{n-1}(\omega)
=\sum_{n'}\chi_{J;n',n}(\omega)\chi^{(k+l)}_{n'-1}(\omega),
\ena
where $\chi_{L}^{(k)}(\omega)\ (L=J,n-1,n'-1)$ are the irreducible characters
of the affine Lie algebra $\widehat{sl_2}$ labeled by the level $k$ and spin
$L/2$.

The one point local height 
probability of the $k-$fusion RSOS model in the regime III 
is given by the same branching 
coefficient as \refeq{charact}. Hence one may regard the space 
${\cal H}_{J;n',n}$ as the  space of states of the $k-$fusion RSOS model
in the regime III, on which Baxter's CTM acts\cite{Bax82,DJKMO}. 

\setcounter{equation}{0}
\section{Discussion}

In summary, we have introduced the elliptic algebra $\uxpslt$.
Based on this algebra, we have extended the bosonization of the ABF model
in \cite{LukPug} to the $k-$fusion RSOS model. We have obtained it as the 
$q-$deformation of the coset conformal field theory $SU(2)_k\times 
SU(2)_{r-k-2}/SU(2)_{r-2}$. A full set of screening currents 
and the highest component of the 
two types of vertex operators have been derived. 
We have observed that these operators give a proper characterization of the 
Fock spaces as the space of states in the $k-$fusion RSOS model. 

We have conjectured that there exists 
a corresponding $q-$deformation of the extended Virasoro algebra
in such a way that 
its screening currents satisfy the algebra $\uxpslt$ and 
the $q-$deformed primary fields are 
determined as the intertwiners between  the infinite representation space of 
$\uxpslt$. In order to establish this point, we need to clarify the
following two points.
\begin{enumerate}
\item
The realization of the $q-$Virasoro generator with the central charge 
\refeq{centralc} and the extra generators. 
\item 
The Hopf algebra  structure of $\uxpslt$. 
\end{enumerate}

Recently, Jimbo has succeeded to derive 
the algebra $\uxpslt$ by the Gauss decomposition of a  central 
extended dynamical $RLL-$relations with the $R-$matrix introduced by 
Enriquez and Felder\cite{EnriFeld}. This result allows us to
clarify 
a Hopf algebra structure of $\uxpslt$.
The work along this line is now in progress.

For the completion of the  identification  of our bosonization with
the $k-$fusion RSOS model, we need further to show that the  
commutation relation among the type I
and the type II vertices are given precisely by the $k-$fused RSOS 
Boltzmann 
weights\cite{FJMMN}. 
We have checked this point both in the  $k=1$ case and in the
case of arbitrary $k$ but for the highest components of the 
type I and type II vertex operators
(Prop. 4.5).  
 
As the scaling limit, we also have obtained a bosonization of the 
algebra $\aheslt$ at arbitrary level $k$. The level one case has shown  
to be a relevant symmetry for deriving the form factors 
in the sine-Gordon theory\cite{Konno97}. In this respect, one may 
relate the higher level case, especially, the case $k=2$ to the 
super sine-Gordon theory.  Our bosonization should be useful for
deriving the form factors of such theory.

It is also an interesting problem to extend our results to the 
higher rank case as well as  to the case of 
other types of Lie algebras. For these 
extensions, the corresponding SOS models are known\cite{JMO,Gepner}.
We expect that these SOS models could be bosonized based on 
the corresponding extension of our $\uxpslt$.

\setcounter{equation}{0}

\section{Acknowledgments} 
The author would like to thank Michio Jimbo for 
stimulating discussions and various suggestions. He is also grateful to
Kenji Iohara, Atsuo Kuniba, Tetsuji Miwa, Stanislav 
Pakuliak and Junichi Shiraishi for valuable discussions.
This work is supported in part by the Ministry of Education Contract 
No.09740028.
\pagebreak[3]

\appendix
\setcounter{equation}{0}
\section{The Algebra ${\cal A}_{\hbar,\eta}(\widehat{sl_2})$}
We here give a brief review of the results in \cite{KLP}.

\begin{dfn}
{The algebra $\aheslt$ is generated by the symbols $\eh, \fh, \th,\ 
\lambda\in\R$
and the central element $c$ with the following relations.}
\bea
&&H^{+}(a)H^{-}(\b)=\frac{\sh\pi\eta({\a-\b-i\hbar(1-\frac{c}{2})})}{
\sh\pi\eta({\a-\b+i\hbar(1+\frac{c}{2})})}
\frac{\sh\pi\eta'({\a-\b+i\hbar(1-\frac{c}{2})})}{
\sh\pi\eta'({\a-\b-i\hbar(1+\frac{c}{2})})}H^{-}(\b)H^{+}(\a),\nonumber\\
\\
&&H^{\pm}(\a)H^{\pm}(\b)=\frac{\sh\pi\eta({\a-\b-i\hbar})}{
\sh\pi\eta({\a-\b+i\hbar})}
\frac{\sh\pi\eta'({\a-\b+i\hbar})}{
\sh\pi\eta'({\a-\b-i\hbar})}
H^{\pm}(\b)H^{\pm}(\a),\\
&&H^{\pm}(\a)E(\b)=\frac{\sh\pi\eta({\a-\b-i\hbar(1\mp\frac{c}{2})})}{
\sh\pi\eta({\a-\b+i\hbar(1\pm\frac{c}{2})})}
E(\b)H^{\pm}(\a),\\
&&H^{\pm}(\a)F(\b)=\frac{\sh\pi\eta'({\a-\b+i\hbar(1\mp\frac{c}{2})})}{
\sh\pi\eta'({\a-\b-i\hbar(1\pm\frac{c}{2})})}
F(\b)H^{\pm}(\a),\\
&&E(\a)E(\b)=\frac{\sh\pi\eta({\a-\b-i\hbar})}
{\sh\pi\eta({\a-\b+i\hbar})}
E(\a)E(\b),\\
&&F(\a)F(\b)=\frac{\sh\pi\eta'({\a-\b+i\hbar})}{
\sh\pi\eta'({\a-\b-i\hbar})}F(\b)F(\a),\\
&&E(\a)F(\b)=F(\b)E(\a)
\nonumber \\
&&\quad =\frac{2\pi}{\hbar}
\Bigl(\frac{1}{\a-\b-\frac{i\hbar c}{2}}H^+\Bigl(\a-\frac{i\hbar c}{4}\Bigr)
-\frac{1}{\a-\b+\frac{i\hbar c}{2}}H^-\Bigl(\b-\frac{i\hbar c}{4}\Bigr)
\Bigr)+O(1),
\ena
where $\hbar$ and $\eta(>0)$ are real  parameter and $1/\eta'-1/\eta=\hbar c>0$. The generating functions (currents)$E(\a), F(\a)$ and $H^\pm(\a)$ are  defined by the following formulae.
\bea
&&E(\a)=\int_{-\infty}^\infty d\lambda e^{i\lambda\a}\eh,\\
&&F(\a)=\int_{-\infty}^\infty d\lambda e^{i\lambda\a}\fh,\\
&&H^\pm(\a)=-\frac{\hbar}{2}\int_{-\infty}^\infty d\lambda e^{i\lambda\a}\th e^{\mp\lambda/2\eta''}
\ena
with $\eta''=\frac{2\eta\eta'}{\eta+\eta'}$.
\end{dfn}

Define another generating functions $e^\pm(\a), f^{\pm}(\a)$ and 
$h^\pm(\a)$ by
\bea
&&e^\pm(\a)={\rm sin}\pi\eta\hbar\int_C\frac{d\gamma}{2\pi i}
\frac{E(\a)}{\sh\pi\eta(\a-\gamma\pm i\hbar\frac{c}{4})},\\
&&f^\pm(\a)={\rm sin}\pi\eta'\hbar\int_{C'}\frac{d\gamma}{2\pi i}
\frac{F(\a)}{\sh\pi\eta'(\a-\gamma\mp i\hbar\frac{c}{4})},\\
&&h^\pm(\a)=\frac{{\rm sin}\pi\eta\hbar}{\pi\eta\hbar}H^{\pm}(\a).
\ena
The following relations hold.
\bea
e^-(\a)=-e^+(\a-i/\eta''),\quad  f^-(\a)=-f^+(\a-i/\eta''),\quad
 h^-(\a)=h^+(\a-i/\eta'').
\ena

The comultiplication of the algebra $\aheslt$ is given 
for $e^+(\a,\xi)=e^+(\a), f^+(\a,\xi)=f^+(\a), h^+(\a,\xi)=h^+(\a)$ 
by the formulae
\bea
&&\Delta c=c'+c''=c\otimes 1 + 1\otimes c, \label{comultic}\\
&&\Delta e^+(\a,\xi)=e^+(\a',\xi)\otimes 1\nonumber\\
&&\qquad\qquad+\sum_{p=0}^\infty
(-1)^p(f^+(\a'-i\hbar,\xi'))^p h^+(\a',\xi')\otimes (e^+(\a'',\xi''))^{p+1},\label{comultie}\\
&&\Delta f^+(\a,\xi)=1\otimes f^+(\a',\xi)\nonumber\\
&&\qquad\qquad
+\sum_{p=0}^\infty
(-1)^p(f^+(\a',\xi'))^{p+1}\otimes \tilde{h}^+(\a'',\xi'') 
(e^+(\a''-i\hbar,\xi''))^{p},\label{comultif}\\
&&\Delta h^+(\a,\xi)
= h^+(\a',\xi')\otimes {h}^+(\a'',\xi'') \nonumber\\
&&\qquad\qquad+
\sum_{p=0}^\infty
(-1)^p[p+1]_\eta(f^+(\a'-i\hbar,\xi'))^p h^+(\a',\xi')\otimes 
{h}^+(\a'',\xi'') (e^+(\a''-i\hbar,\xi''))^{p},\nonumber\\
\label{comultih}
\ena
where $\xi=1/\eta=\xi',\ \xi''=\xi+\hbar c$, $\a'=\a+i\hbar c''/4,\ 
\a''=\a-i\hbar c'/4$ and 
$$
[p]_\eta=\frac{{\rm sin}\pi\eta\hbar p}{{\rm sin}\pi\eta\hbar},
\qquad
\tilde{h}^+(\a)=\frac{\eta}{\eta'}
\frac{{\rm sin}\pi\eta'\hbar }{{\rm sin}\pi\eta\hbar} h^+(\a).
$$

Let $V^{(l)}$ be the $l+1$-dimensional  representation of 
$U_\rho(sl_2), \ \rho=e^{i\pi\eta\hbar}$ with a basis 
$v_m\ (m=0,1,2,..,l)$. The $l+1$-dimensional evaluation representation 
$V^{(l)}(\z)=V^{(l)}\otimes \C[[e^{i\lambda\z}]],\ \lambda\in\R,\ \z\in\C$ 
is given for the base
$v_m(\z)\in V^{(l)}(\z)$ by 
\bea
&&e^+(\a)v_m(\z)=-\frac{\sh i\pi\eta\hbar m}
{\sh\pi\eta(\a-\z+i\hbar \frac{l-2m-1}{2})}v_{m-1}(\z),\label{evale}\\
&&f^+(\a)v_m(\z)=-\frac{\sh\pi\eta\hbar(l- m)}
{\sh\pi\eta(\a-\z+i\hbar \frac{l-2m-1}{2})}v_{m+1}(\z),\label{evalf}\\
&&h^+(\a)v_m(\z)=\frac{\sh\pi\eta(\a-\z-i\hbar\frac{l+1}{2})\sh\pi\eta(\a-\z+i\hbar\frac{l+1}{2})}
{\sh\pi\eta(\a-\z+i\hbar \frac{l-2m+1}{2})
\sh\pi\eta(\a-\z+i\hbar \frac{l-2m-1}{2})}v_{m}(\z).\label{evalh}
\ena

An infinite dimensional highest weight representation of the algebra 
$\aheslt$ is constructed as a Fock space $\F$ of free bosons\cite{JKM,KLP}. 
See also Sec.6. Here the highest weight property means that
\bea
&&\eh\ket{h.w.}=0,\quad \fh\ket{h.w.}=0 \qquad \lambda\in\R.
\ena

There are two types of intertwining operators (vertex operators):
\bea
&&{\rm Type\quad I}\quad \Phi^{(l)}(\z)\ :\ \F\to\F\otimes V^{(l)}(\z+
\frac{i\hbar}{2}),\label{typei}\\
&&{\rm Type\quad II}\quad \Psi^{(l)*}(\z)\ :\ V^{(l)}(\z+\frac{i\hbar}{2})\otimes\F\to\F \label{typeii} 
\ena
satisfying
\bea
&&\Phi^{(l)}(\z)\ x=\Delta(x)\Phi^{(l)}(\z),\\
&&\Psi^{(l)*}(\z)\Delta(x)=\ x \Psi^{(l)*}(\z).
\ena

The components of the intertwiners are defined as follows.
\bea
&&\Phi^{(l)}(\z)u=\sum_{m=0}^l\Phi^{(l)}_m(\z)u\otimes v_m ,\\
&&\Psi^{(l)*}(\z)(v_m\otimes u)=\Psi^{(l)}_m(\z)u\quad u\in\F.
\ena

Using the comultiplication \refeq{comultic}--\refeq{comultih} and the finite dimensional evaluation 
representation \refeq{evale}--\refeq{evalh}, the intertwining relations are rewritten as follows.
For type I,
\bea
&&\Phil{l}h^{\pm}(\a)=\frac{\sh\pi\eta'(\a-\z+\frac{i\hbar}{2}(l\mp 
\frac{k}{2}))}
{\sh\pi\eta'(\a-\z-\frac{i\hbar}{2}(l\pm \frac{k}{2}))}h^{\pm}(\a)\Phil{l},
\label{inthi}\\
&&[e^{\pm}(\a),\Phil{l}]=0,\label{intei}\\
&&\sh i\pi \eta\hbar(l-m+1)\Phil{m-1}\nonumber\\
&&\quad =-\sh\pi\eta'(\a-\z+\frac{i\hbar}{2}(l-2m+2\mp \frac{k}{2}))
\Bigl[\Phil{m}f^\pm(\a)\nonumber\\
&& -
\frac{\sh\pi\eta'(\a-\z-\frac{i\hbar}{2}(l+2\pm \frac{k}{2}))\sh\pi\eta'(\a-\z+i\hbar(l\mp 
k/2))}
{\sh\pi\eta'(\a-\z+\frac{i\hbar}{2}(l-2m\mp \frac{k}{2}))\sh\pi\eta'(\a-\z+i\hbar(l-2m-2\mp
\frac{k}{2}))}
f^\pm(\a)\Phil{m}\Bigr],\nonumber \\
\label{intfi}
\ena
and for type II,
\bea
&&h^{\pm}(\a)\Psils{l}=\frac{\sh\pi\eta(\a-\z+\frac{i\hbar}{2}(l\pm \frac{k}{2}))}
{\sh\pi\eta(\a-\z-\frac{i\hbar}{2}(l\mp \frac{k}{2}))}\Psils{l}h^{\pm}(\a),
\label{inthii}\\
&&[f^{\pm}(\a),\Psils{l}]=0,\label{intfii}\\
&&\sh i\pi \eta\hbar m\Psils{m-1}\nonumber\\
&&\quad =-\sh\pi\eta(\a-\z+\frac{i\hbar}{2}(l-2m-2\pm \frac{k}{2}))
e^\pm(\a)\Psils{m}\nonumber\\
&&+\frac{\sh\pi\eta(\a-\z-\frac{i\hbar}{2}(l+2\mp \frac{k}{2}))
\sh\pi\eta(\a-\z+\frac{i\hbar}{2}(l\pm \frac{k}{2}))}
{\sh\pi\eta(\a-\z+\frac{i\hbar}{2}(l-2m\pm \frac{k}{2}))}
\Psils{m}e^\pm(\a)\label{inteii},
\ena

In some unknown reason, the intertwiners relevant for the physical problems 
such as the mass less XXZ model\cite{JKM} and the sine-Gordon 
theory\cite{Lukyanov,Konno97} prefer the twisted intertwining relations. 
We denote by $\tilde\Phi^{(l)}(\z)$ and $\tilde\Psi^{(l)*}(\z)$  the 
intertwiners of the same types as \refeq{typei} and \refeq{typeii} but
obeying the following twisted intertwining relations.
\bea
&&\tilde\Phi^{(l)}(\z)\ \iota(x)=\Delta(x)\tilde\Phi^{(l)}(\z),\\
&&\tilde\Psi^{(l)*}(\z)\Delta(x)=\ \iota(x) \tilde\Psi^{(l)*}(\z),\\
\ena
where $\iota$ is the following involution of $\aheslt$
\bea
&&\iota(\eh)=-\eh,\quad \iota(\fh)=-\fh,\quad  \iota(\th)=\th.
\ena
Then the relations \refeq{inthi} and \refeq{inthii} remains the same
but \refeq{intei}, \refeq{intfi}, \refeq{intfii} and \refeq{inteii} are
replaced with
\bea
&&\{e^{\pm}(\a),\tilde\Phil{l}\}=0,\label{tintei}\\
&&\sh i\pi \eta\hbar(l-m+1)\tilde\Phil{m-1}\nonumber\\
&&\quad =\sh\pi\eta'(\a-\z+\frac{i\hbar}{2}(l-2m+2\mp \frac{k}{2}))
\Bigl[\tilde\Phil{m}f^\pm(\a)\nonumber\\
&& +
\frac{\sh\pi\eta'(\a-\z-\frac{i\hbar}{2}(l+2\pm \frac{k}{2}))\sh\pi\eta'(\a-\z+i\hbar(l\mp 
k/2))}
{\sh\pi\eta'(\a-\z+\frac{i\hbar}{2}(l-2m\mp \frac{k}{2}))\sh\pi\eta'(\a-\z+i\hbar(l-2m-2\mp
\frac{k}{2}))}
f^\pm(\a)\tilde\Phil{m}\Bigr],\nonumber\\
\label{tintfi} \\
&&\{f^{\pm}(\a),\tilde\Psils{l}\}=0,\label{tintfii}\\
&&\sh i\pi \eta\hbar m\tilde\Psils{m-1}\nonumber\\
&&\quad =\sh\pi\eta(\a-\z+\frac{i\hbar}{2}(l-2m-2\pm \frac{k}{2}))
e^\pm(\a)\tilde\Psils{m}\nonumber\\
&&\qquad+\frac{\sh\pi\eta(\a-\z-\frac{i\hbar}{2}(l+2\mp \frac{k}{2}))
\sh\pi\eta(\a-\z+\frac{i\hbar}{2}(l\pm \frac{k}{2}))}
{\sh\pi\eta(\a-\z+\frac{i\hbar}{2}(l-2m\pm \frac{k}{2}))}
\tilde\Psils{m}e^\pm(\a),\label{tinteii} 
\ena

\end{document}